\documentclass[11pt]{article}
\usepackage{amssymb}
\usepackage{amsmath}
\usepackage{t1enc}
\usepackage[latin1]{inputenc}
\usepackage[english]{babel}
\usepackage{amsmath,amsthm}
\usepackage{amsfonts}
\usepackage{latexsym}
\usepackage[dvips]{graphicx}
\usepackage{graphicx}
\usepackage[natural]{xcolor}
\usepackage{setspace}
\usepackage{fullpage}

\newtheorem{definition}{Definition}
\newtheorem{theorem}{Theorem}
\newtheorem{lemma}{Lemma}

\newtheorem{example}{Example}
\newtheorem{remark}{Remark}

\makeatletter
   
   \@addtoreset{equation}{section}
\makeatother

\begin{document}

\title{Necessary Condition for Near Optimal Control of Linear
Forward-backward Stochastic Differential Equations}
\author{Liangquan Zhang$^{1}$\thanks{
L. Zhang acknowledges the financial support partly by the National Nature
Science Foundation of China (No. 11201263 and No. 11201264) and the Nature
Science Foundation of Shandong Province (No. ZR2012AQ004). E-mail:
xiaoquan51011@163.com.}, \quad Jianhui Huang$^{2}$\thanks{
J. Huang acknowledges the financial support partly by Hong Kong RGC grant
501010. E-mail: majhuang@polyu.edu.hk.}, \quad Xun Li$^{2}$ \thanks{
X. Li acknowledges the financial support partly by Hong Kong RGC grant
524109. E-mail: malixun@polyu.edu.hk.} \\
{\small 1. INRIA, Campus de Beaulieu, 35042 Rennes, France.}\\
{\small 2. Department of Applied Mathematics} \\
{\small \ \ The Hong Kong Polytechnic University, Hung Hom, Kowloon, Hong
Kong, China.}}
\maketitle

\begin{abstract}
This paper investigates the near optimal control for a kind of linear
stochastic control systems governed by the forward backward stochastic
differential equations, where both the drift and diffusion terms are allowed
to depend on controls and the control domain is not assumed to be convex. In
the previous work (Theorem 3.1) of the second and third authors [\textit{%
Automatica} \textbf{46} (2010) 397-404], some problem of near optimal
control with the control dependent diffusion is addressed and our current
paper can be viewed as some direct response to it. The necessary condition
of the near-optimality is established within the framework of optimality
variational principle developed by Yong [\textit{SIAM J. Control Optim.}
\textbf{48} (2010) 4119--4156] and obtained by the convergence technique to
treat the optimal control of FBSDEs in unbounded control domains by Wu [%
\textit{Automatica} \textbf{49} (2013) 1473--1480]. Some new estimates are
given here to handle the near optimality. In addition, an illustrating
example is discussed as well.
\end{abstract}



\section{Introduction}

\label{sect:1}

Due to the nature of uncertainty, solutions to a forward stochastic system
governed by Itô-based \textit{stochastic differential equations} (SDEs in
short) need to be non-anticipative. The equation for a conventional Itô SDEs
can be naturally solved in a forward-looking way by starting with the
initial state. In some financial engineering problems, however, it is
inherent that some terminal states are specified and one must consider a
stochastic dynamic system in a backward fashion. For example, one needs to
determine the option price with a given terminal payoff (as a random
variable on the underlying asset). This results in a \textit{backward
stochastic differential equations} (BSDEs in short) with a terminal
condition. This theory can be traced back to Bismut \cite{B} who studied
linear BSDEs motivated by stochastic control problems, and Pardoux and Peng
\cite{PP} who proved the well-posedness for nonlinear BSDEs. Since then,
BSDEs have been extensively studied and used in the areas of applied
probability and optimal stochastic controls, particularly in financial
engineering. Moreover, initiated by Antonelli \cite{A}, \textit{%
forward-backward stochastic differential equations} (FBSDEs in short) have
also been investigated systematically. For instance, Ma, Protter and Yong
\cite{MPY} established the four-step-scheme, Pardoux and Tang \cite{PT}, Hu
and Peng \cite{HP} , Peng and Wu \cite{PW} and Yong \cite{Y1997} developed
the method of continuation, Huang, Li and Wang \cite{HLW,HLW1} analyzed
partial information control problems using FBSDEs, Lim and Zhou \cite{LZ}
formulated and solved backward linear-quadratic controls, and Yong \cite%
{Y,Y2010} further considered coupled FBSDEs with mixed initial-terminal
conditions. See also relevant work by Wu \cite{W, W1, WY}, more references
therein.


Near optimization has been investigated by many literatures for both theory
and applications. On the one hand, near optimal controls are more available
than optimal ones. Indeed, optimal controls may not exist in lots of
situations, while near optimal controls always exist, and it is much easier
to derive the near optimal controls than optimal ones, both analytically and
numerically. On the other hand, since the there are many candidates for near
optimal controls, it is possible to select among them appropriates ones that
are easier for analysis and implementation (see \cite{ZS} reference therein).

As a matter of fact, the near optimal control for the forward deterministic
and stochastic systems have been extensively studied. We refer the reader to
the monographs \cite{CM, Y, Z1995, Z1996, Z1998, ZS} for deterministic and
stochastic cases. Bahlali, Khelfallah and Mezerdi \cite{BKB} investigated
the near optimal control of FBSDEs (see also references therein Hafayed et
al. \cite{HA1, HAA, HA2, HAV, HVA, HVA1}). Based on Ekeland's principle and
spike variation, a necessary and sufficient condition of near optimality for
the near optimal control are established. However, in their work, the
diffusion coefficient is independent of the control variable. The similar
hypothesis was put in the work of Huang, Li and Wang \cite{HLW} for linear
case. Besides, Hui, Huang, Li and Wang \cite{HHLW} also considered the near
optimal control for general form of FBSDEs, with the assumption that the
control domain is convex. It is remarkable that some problem to near optimal
control is addressed in previous work of second and third authors (see \cite%
{HLW}) in which the diffusion term depends on the control. The difficulty to
this problem when the controlled systems are FBSDEs is also discussed. Our
aims in this paper is to fulfill this research gap by removing this
assumption, that is, \textit{the diffusion coefficient is independent of the
control variable and control domain is convex}. Our methods are mainly based
on the Ekeland's principle, spike variation and reduction technique
developed recently by Yong \cite{Y} and the methodology recently introduced
by Wu \cite{W} to consider the optimal control problem for FBSDEs in the
general case of control domains including Lagrange multipliers.

Let us make it more precise. First of all, we introduce the controlled
initial value problem for a system of SDEs, where the pair $\left( x\left(
\cdot \right) ,y\left( \cdot \right) \right) $ is regarded as the state
process and $\left( z\left( \cdot \right) ,u\left( \cdot \right) \right) $
is regarded as the control process in bounded control domains. Meanwhile, we
regard the original terminal condition $y\left( T\right) =Mx\left( T\right) $
as the terminal state constraint. Next it is possible to translate the near
optimal control Problem $(\tilde{C}^{\varepsilon })$ into a high-dimensional
reduced near optimal control problem driven by the standard SDEs with state
constraint (for more information see Problem $(C^{\varepsilon })$ in Section %
\ref{sect:3}). We mention that the advantage of this reduced near optimal
control problem is that one needs not much regularity/integrability of
process $z\left( \cdot \right) $ since it is treated as a control process.
Hence, it is possible to apply the Ekeland's variational principle to handle
this high-dimensional reduced near optimal control problem with state
constraint. Afterwards, the necessary conditions for the near optimal
control of Problem $(C)$ are derived by Problem $(\tilde{C}^{\varepsilon })$%
. Finally, by convergence technique we obtain the general case of control
domains and complete our proof.

The paper is organized as follows. The notations, preliminaries and some
basic definitions are given in Section \ref{sect:2}. In Section \ref{sect:3}%
, under some suitable assumptions, we state the main result of this paper,
together with some discussions of special cases. The application of our
theoretical results will be shown in Section \ref{sect:4}. Some conclusion
is given in Section \ref{sect:5}. Finally, we present some technique proofs
in Appendix. For the simplicity of notations, we consider the case where
both $x$ and $y$ are one-dimensional, and the control $u$ is also
one-dimensional.

\section{Notation and preliminaries}

\label{sect:2}

Throughout this paper, we denoted by $\mathbb{R}$ the space of
one-dimensional Euclidean space, by $\mathbb{R}^{n\times d}$ the space the
matrices with order $n\times d$, by $\mathbf{S}^{n}$ the space of symmetric
matrices with order $n\times n$. $\left\langle \cdot ,\cdot \right\rangle $
and $\left\vert \cdot \right\vert $ denote the scalar product and norm in
the Euclidean space, respectively. $\top $ is the transpose of a matrix.

Let $\mathbb{U}$ be a given set in some Euclidean space $\mathbb{R}$. Let $%
(\Omega ,\mathcal{F},\{\mathcal{F}_{t}\}_{t\geq 0},P)$ be a complete
filtered probability space on which a one-dimensional standard Brownian
motion $W(\cdot )$ is defined, with $\{\mathcal{F}_{t}\}_{t\geq 0}$ being
its natural filtration, augmented by all the $P$-null sets.

We now introduce the following spaces of process:
\begin{align*}
\mathcal{S}^{2}(0,1;\mathbb{R})\triangleq & \left\{ \mathbb{R}\text{-valued }%
\mathcal{F}_{t}\text{-adapted process }\phi (t)\text{; }\mathbb{E}\left[
\sup\limits_{0\leq t\leq 1}\left\vert \phi (t)\right\vert ^{2}\right]
<\infty \right\} , \\
\mathcal{M}^{2}(0,1;\mathbb{R})\triangleq & \left\{ \mathbb{R}\text{-valued }%
\mathcal{F}_{t}\text{-adapted process }\varphi (t)\text{; }\mathbb{E}\left[
\int_{0}^{1}\left\vert \varphi (t)\right\vert ^{2}\mbox{\rm d}t\right]
<\infty \right\} ,
\end{align*}%
and denote $\mathcal{N}^{2}\left[ 0,1\right] =\mathcal{S}^{2}(0,1;\mathbb{R}%
)\times \mathcal{S}^{2}(0,1;\mathbb{R})\times \mathcal{M}^{2}(0,1;\mathbb{R}%
).$ Clearly, $\mathcal{N}^{2}\left[ 0,1\right] $ forms a Banach space. Any
process in $\mathcal{N}^{2}\left[ 0,1\right] $ is denoted by $\Theta (\cdot
)=(x(\cdot ),y(\cdot ),z(\cdot ))$, whose norm is given by
\begin{equation*}
\left\Vert \Theta (\cdot )\right\Vert _{\mathcal{N}^{2}\left[ 0,1\right] }=%
\mathbb{E}\left[ \sup\limits_{t\in \left[ 0,1\right] }\left\vert
x(t)\right\vert ^{2}+\sup\limits_{t\in \left[ 0,1\right] }\left\vert
y(t)\right\vert ^{2}+\int_{0}^{1}\left\vert z(t)\right\vert ^{2}\mbox{\rm d}t%
\right] .
\end{equation*}

\subsection{Formulation of Near Optimal Control Problem and Basic Assumptions%
}

We study the stochastic control systems which are described by a linear
FBSDEs of the type:
\begin{equation}
\left\{
\begin{array}{rcl}
\mbox{\rm d}x(t) & = & \left[ A(t)x(t)+B(t)u(t)\right] \mbox{\rm d}t+\left[
C(t)x(t)+D(t)u(t)\right] \mbox{\rm d}W(t), \\[2mm]
\mbox{\rm d}y(t) & = & -\left[ a(t)x(t)+b(t)y(t)+c(t)u(t)\right] \mbox{\rm d}%
t+z(t)\mbox{\rm d}W(t), \\[2mm]
x(0) & = & x_{0},\qquad y(1)=Mx(1),%
\end{array}%
\right.  \label{2.1}
\end{equation}%
where $A(\cdot ),B(\cdot ),C(\cdot ),D(\cdot ),a(\cdot ),b(\cdot )$ and $%
c(\cdot )$ are bounded deterministic functions with values in $\mathbb{R}$, $%
M$ is a constant, and $u(\cdot )$ is a control process.

The control process $u(\cdot ):\left[ 0,1\right] \times \Omega \rightarrow
\mathbb{U}$ is called admissible, if it is an $\mathcal{F}_{t}$-adapted
process with values in $\mathbb{U}$. The set of all admissible controls is
denoted by $\mathcal{U}_{ad}[0,1]$.

Under the above assumptions, for any $u(\cdot )\in \mathcal{U}_{ad}\left[ 0,1%
\right] $, it is easy to check that FBSDEs (\ref{2.1}) admit a unique $%
\mathcal{F}_{t}$-adapted solution denoted by the triple $(x(\cdot ),y(\cdot
),z(\cdot ))\in \mathcal{S}^{2}(0,1;\mathbb{R})\times \mathcal{S}^{2}(0,1;%
\mathbb{R})\times \mathcal{M}^{2}(0,1;\mathbb{R})$.

\noindent The cost functional is given by
\begin{equation}
J(u(\cdot ))=\mathbb{E}\left[ \int_{0}^{1}l(t,x(t),y(t),u(t))\mbox{\rm d}%
t+\phi (x(1))+\gamma (y(0))\right] ,  \label{2.2}
\end{equation}%
where
\begin{equation*}
\begin{array}{rl}
\phi : & \mathbb{R}\rightarrow \mathbb{R}, \\
\gamma : & \mathbb{R}\rightarrow \mathbb{R}, \\
l: & \left[ 0,1\right] \times \mathbb{R\times R}\times \mathbb{U}\rightarrow
\mathbb{R}.%
\end{array}%
\end{equation*}%
The classical object of the optimal control problem is to minimize the cost
functional $J(u(\cdot)) $, over all $u(\cdot) \in \mathcal{U}_{ad}\left[ 0,1%
\right]$. We denote the above problem by $(C)$.

\noindent \textbf{Problem $(C)$}. Find $\bar{u}(\cdot )\in \mathcal{U}%
_{ad}[0,1]$, such that
\begin{equation}
J(\bar{u}(\cdot ))=\inf\limits_{u(\cdot )\in \mathcal{U}_{ad}\left[ 0,1%
\right] }J(u(\cdot )).  \label{2.3}
\end{equation}%
Any $\bar{u}(\cdot )\in \mathcal{U}_{ad}\left[ 0,1\right] $ satisfying (\ref%
{2.3}) is called an optimal control process of Problem $(C)$, and the
corresponding state process, denoted by $(\bar{x}(\cdot ),\bar{y}(\cdot ),%
\bar{z}(\cdot ))$, is called optimal state process. We also refer to $(\bar{x%
}(\cdot ),\bar{y}(\cdot ),\bar{z}(\cdot ),\bar{u}(\cdot ))$ as an optimal
4-tuple of Problem $(C)$.

However, the control problem under consideration in this paper is to find
the a control in $\mathcal{U}_{ad}\left[ 0,1\right] $, which minimizes or
``nearly'' minimizes $J( \bar{u}( \cdot ) ) $ over $\mathcal{U}_{ad}[0,1]$.
From this point, we need the following definitions.

\begin{definition}[\textbf{Optimal Control}]
Any admissible control $\bar{u}( \cdot ) \in \mathcal{U}_{ad} \left[ 0,1%
\right] ,$ is called optimal$,$ if $\bar{u}( \cdot ) $ attains the minimum
of $J( u( \cdot ) ) .$
\end{definition}

\begin{definition}[$\protect\varepsilon $\textbf{-Optimal Control}]
For a given $\varepsilon >0,$ an admissible control $u^{\varepsilon }( \cdot
) $ is called $\varepsilon $-optimal if
\begin{equation*}
\left\vert J( u^{\varepsilon }( \cdot ) ) -J( \bar{u}( \cdot ) ) \right\vert
\leq \varepsilon .
\end{equation*}
\end{definition}

\begin{definition}
Both a family of admissible controls $\left\{ u^{\varepsilon }(\cdot
)\right\} $ parameterized by $\varepsilon >0$ and any element $%
u^{\varepsilon }(\cdot ),$ in the family, are called near optimal if
\begin{equation*}
\left\vert J(u^{\varepsilon }(\cdot ))-J(\bar{u}(\cdot ))\right\vert \leq
r(\varepsilon )
\end{equation*}%
holds for sufficient small $\varepsilon ,$ where $r$ is a function of $%
\varepsilon $ satisfying $r(\varepsilon )\rightarrow 0$ as $\varepsilon
\rightarrow 0.$ The estimate $r(\varepsilon )$ is called an error bound. If $%
r(\varepsilon )=C\varepsilon ^{\delta }$ for some $\delta >0$ independent of
the constant $C$, then $u^{\varepsilon }(\cdot )$ is called near optimal
with order $\varepsilon ^{\delta }$.
\end{definition}

\noindent \textbf{Problem $(C^{\varepsilon })$}. Find $\bar{u}^{\varepsilon
}(\cdot )\in \mathcal{U}_{ad}[0,1]$, such that
\begin{equation}
J(\bar{u}(\cdot ))=\inf\limits_{u(\cdot )\in \mathcal{U}_{ad}\left[ 0,1%
\right] }J(u(\cdot ))+\varepsilon .  \label{2.31}
\end{equation}

\noindent Any $\bar{u}^{\varepsilon }(\cdot )\in \mathcal{U}_{ad}\left[ 0,1%
\right] $ satisfying (\ref{2.31}) is called a near optimal control process
of Problem $(C)$, and the corresponding state process, denoted by $(\bar{x}%
^{\varepsilon }(\cdot ),\bar{y}^{\varepsilon }(\cdot ),\bar{z}^{\varepsilon
}(\cdot ))$, is called optimal state process. We also refer to $(\bar{x}%
^{\varepsilon }(\cdot ),\bar{y}^{\varepsilon }(\cdot ),\bar{z}^{\varepsilon
}(\cdot ),\bar{u}^{\varepsilon }(\cdot ))$ as an optimal 4-tuple of Problem $%
(C)$.

Hereafter, $C>0$ stands for a generic constant which can be different at
different places.

\section{Main Result}

\label{sect:3}

\subsection{Necessary Condition of Near Optimality}

In this section, we first present our necessary conditions for the near
optimal control of Problem $(C)$ under some suitable assumptions. Due to the
assumptions introduced in Section \ref{sect:2}. There exists a constant $L>0$
such that%
\begin{align*}
\left\vert A(t)(x-x^{\prime })\right\vert ^{2}+\left\vert C(t)(x-x^{\prime
})\right\vert ^{2}+\left\vert a(t)(x-x^{\prime })+b(t)(y-y^{\prime
})\right\vert ^{2}& \leq L(\left\vert x-x^{\prime }\right\vert
^{2}+\left\vert y-y^{\prime }\right\vert ^{2}), \\
\forall t& \in \left[ 0,1\right] ,\text{ }(x,y),\text{ }(x^{\prime
},y^{\prime })\in \mathbb{R}\times \mathbb{R},
\end{align*}%
and
\begin{align*}
\left\vert A(t)x+B(t)u\right\vert ^{2}+\left\vert C(t)x+D(t)u\right\vert
^{2}+\left\vert a(t)x+b(t)y+c(t)u\right\vert ^{2}& \leq L(1+\left\vert
x\right\vert ^{2}+\left\vert y\right\vert ^{2}), \\
\forall (t,z,u)& \in \left[ 0,1\right] \times \mathbb{R}\times \mathbb{U},%
\text{ }(x,y)\in \mathbb{R}\times \mathbb{R}.
\end{align*}

\noindent To establish the necessary condition, we need the following
assumption:

\begin{enumerate}
\item[\textbf{(H1)}] The maps $\phi ,$ $\gamma $ are twice continuously
differentiable with respect to $\left( x,y\right) $. $l_{x},$ $l_{y},$ $\phi
_{x}$ and $\gamma _{y}$ grow linearly about $\left( x,y,u\right) $ and is
continuous in $\left( t,u\right) .$ Moreover, $l_{xx},$ $l_{yy},$ $l_{xy},$ $%
\phi _{xx}$ and $\gamma _{yy}$ are bounded.
\end{enumerate}

\noindent Now, let $(x^{\varepsilon }(\cdot ),y^{\varepsilon }(\cdot
),z^{\varepsilon }(\cdot ),u^{\varepsilon }(\cdot ))$ be a near optimal
4-tuple of Problem $(C^{\varepsilon })$. We introduce
\begin{align*}
\mathcal{B}_{X}(t,\cdot )& \triangleq \left(
\begin{array}{cc}
A(t) & 0 \\
-a(t) & -b(t)%
\end{array}%
\right) , \\
\Sigma _{X}(t,\cdot )& \triangleq \left(
\begin{array}{cc}
C(t) & 0 \\
0 & 0%
\end{array}%
\right) .
\end{align*}%
Our main result of this paper is following:

\begin{theorem}
\label{thm:4} Suppose (H1) holds. Then, for any $\beta \in \lbrack 0,\frac{1%
}{3})$, there exist a constant $C_{1}=C_{1}(\beta )$ such that for any fixed
$\varepsilon >0$ and any $\varepsilon $-optimal $(x^{\varepsilon }(\cdot
),y^{\varepsilon }(\cdot ),z^{\varepsilon }(\cdot ),u^{\varepsilon }(\cdot
)) $ of the problem $(C)$, there exist two parameters $\theta
_{0}^{\varepsilon }$ and $\theta _{1}^{\varepsilon }$ ($\mathcal{F}_{1}$%
-measurable random variable) with $\left\vert \theta _{0}^{\varepsilon
}\right\vert ^{2}+\mathbb{E}\left\vert \theta _{1}^{\varepsilon }\right\vert
^{2}=1$, $\theta _{0}^{\varepsilon }\geq 0$ holds that
\begin{align}
& \int_{0}^{1}\bigg[\left\langle p^{\varepsilon }(t),B(t)(u-u^{\varepsilon
}(t))\right\rangle +\left\langle k^{\varepsilon }(t),D(t)(u-u^{\varepsilon
}(t))\right\rangle -\left\langle q^{\varepsilon }(t),c(t)(u-u^{\varepsilon
}(t))\right\rangle  \notag \\
& +\theta _{0}^{\varepsilon }\left[ l(t,x^{\varepsilon }(t),y^{\varepsilon
}(t),u)-l(t,x^{\varepsilon }(t),y^{\varepsilon }(t),u^{\varepsilon }(t))%
\right] +\frac{1}{2}D^{\varepsilon }(t)(u-u^{\varepsilon }(t))^{2}P_{1}(t)%
\bigg]\mbox{\rm d}t\geq -C_{1}\theta _{0}^{\varepsilon }\varepsilon ^{\beta
},  \label{3.3}
\end{align}%
where
\begin{equation}
\left\{
\begin{array}{rcl}
-\mbox{\rm d}p^{\varepsilon }(t) & = & \left[ A(t)p^{\varepsilon
}(t)-a(t)q^{\varepsilon }(t)+C(t)k^{\varepsilon }(t)+\theta
_{0}^{\varepsilon }l_{x}^{\varepsilon }(t,\cdot )\right] \text{d}%
t-k^{\varepsilon }(t)\mbox{\rm d}W(t), \\[2mm]
\text{ \ }\mbox{\rm d}q^{\varepsilon }(t) & = & \left[ -b(t)q^{\varepsilon
}(t)-\theta _{0}^{\varepsilon }l_{y}^{\varepsilon }(t,\cdot )\right] %
\mbox{\rm d}t, \\[2mm]
\text{ \ }p^{\varepsilon }(1) & = & \theta _{0}^{\varepsilon }\phi
_{x}(x^{\varepsilon }(1))-M\theta _{1}^{\varepsilon },\qquad q^{\varepsilon
}(0)=-\theta _{0}^{\varepsilon }\mathbb{\gamma }_{y}(y^{\varepsilon }(0)),%
\end{array}%
\right.  \label{3.4}
\end{equation}%
and
\begin{equation}
\left\{
\begin{array}{rcl}
-\mbox{\rm d}P^{\varepsilon }(t) & = & \big[\mathcal{B}_{X}(t,\cdot )^{\top
}P^{\varepsilon }(t)+P^{\varepsilon }(t)\mathcal{B}_{X}(t,\cdot )+\Sigma
_{X}(t,\cdot )^{\top }P^{\varepsilon }(t)\Sigma _{X}(t,\cdot ) \\[2mm]
&  & +\Sigma _{X}(t,\cdot )^{\top }Q^{\varepsilon }(t)+Q^{\varepsilon
}(t)\Sigma _{X}(t,\cdot )+H_{XX}^{\varepsilon }(t,\cdot )\big]\mbox{\rm d}%
t-Q^{\varepsilon }(t)\mbox{\rm d}W(t), \\[2mm]
\text{ \ }P^{\varepsilon }(1) & = & \left(
\begin{array}{cc}
\theta _{0}^{\varepsilon }\phi _{xx}(x^{\varepsilon }(1)) & 0 \\
0 & 0%
\end{array}%
\right) ,%
\end{array}%
\right.  \label{3.5}
\end{equation}%
where%
\begin{equation*}
\left\{
\begin{array}{lll}
\text{ \ \ }l_{x}^{\varepsilon }(t,\cdot ) & = & l_{x}\left(
t,x^{\varepsilon }(t),y^{\varepsilon }(t),u^{\varepsilon }(t)\right) , \\
H_{XX}^{\varepsilon }(t,\cdot ) & = & H_{XX}(t,x^{\varepsilon
}(t),y^{\varepsilon }(t),u^{\varepsilon }(t),p^{\varepsilon
}(t),q^{\varepsilon }(t),k^{\varepsilon }(t),\theta _{0}^{\varepsilon }),%
\end{array}%
\right.
\end{equation*}%
and the Hamiltonian function $H:\left[ 0,T\right] \times \mathbb{R}\times
\mathbb{R}\times \mathbb{U}\times \mathbb{R}\times \mathbb{R}\times \mathbb{%
R\times R\rightarrow R}$ is defined as follows:%
\begin{align*}
H(t,x,y,u,p,q,k,\theta )\triangleq & \left\langle p,A(t)x+B(t)u\right\rangle
-\left\langle q,a(t)x+b(t)y+c(t)u\right\rangle \\
& +\left\langle k,C(t)x+D(t)u\right\rangle +\theta l(t,x,y,u).
\end{align*}
\end{theorem}

The proof can be seen in Appendix. Some remarks are in order.

\begin{remark}
\textsl{Actually, the second order adjoint equations (\ref{3.5}) can be
rewritten as the following three BSDEs if we introduce that
\begin{equation*}
P^{\varepsilon }\left( \cdot \right) \triangleq \left(
\begin{array}{cc}
P_{1}^{\varepsilon }(\cdot ) & P_{2}^{\varepsilon }(\cdot ) \\
P_{2}^{\varepsilon }(\cdot ) & P_{3}^{\varepsilon }(\cdot )%
\end{array}%
\right) ,\text{ }Q^{\varepsilon }(\cdot )\triangleq \left(
\begin{array}{cc}
Q_{1}^{\varepsilon }(\cdot ) & Q_{2}^{\varepsilon }(\cdot ) \\
Q_{2}^{\varepsilon }(\cdot ) & Q_{3}^{\varepsilon }(\cdot )%
\end{array}%
\right) ,
\end{equation*}%
then we have
\begin{equation*}
\left\{
\begin{array}{rcl}
-\mbox{\rm d}P_{1}^{\varepsilon }(t) & = & \left[ 2A(t)P_{1}^{\varepsilon
}(t)+2C(t)P_{1}^{\varepsilon }(t)+2C(t)Q_{1}^{\varepsilon
}(t)-2a(t)P_{2}^{\varepsilon }(t)-\theta _{0}^{\varepsilon
}l_{xx}^{\varepsilon }(t,\cdot )\right] \mbox{\rm d}t-Q_{1}^{\varepsilon }(t)%
\mbox{\rm d}W(t), \\[2mm]
P_{1}^{\varepsilon }(1) & = & \theta _{0}^{\varepsilon }\phi
_{xx}(x^{\varepsilon }(1)),%
\end{array}%
\right.
\end{equation*}%
\begin{equation*}
\left\{
\begin{array}{rcl}
-\mbox{\rm d}P_{2}^{\varepsilon }(t) & = & \left[ A(t)P_{2}^{\varepsilon
}(t)+C(t)Q_{2}^{\varepsilon }(t)-a(t)P_{3}^{\varepsilon
}(t)-P_{2}^{\varepsilon }(t)b(t)-\theta _{0}^{\varepsilon
}l_{xy}^{\varepsilon }(t,\cdot )\right] \mbox{\rm d}t-Q_{2}^{\varepsilon }(t)%
\mbox{\rm d}W(t), \\[2mm]
P_{2}^{\varepsilon }(1) & = & 0,%
\end{array}%
\right.
\end{equation*}%
and
\begin{equation*}
\left\{
\begin{array}{rcl}
-\mbox{\rm d}P_{3}^{\varepsilon }(t) & = & \left[ -2b(t)P_{3}^{\varepsilon
}(t)-\theta _{0}^{\varepsilon }l_{yy}^{\varepsilon }(t,\cdot )\right]
\mbox{\rm
d}t-Q_{3}^{\varepsilon }(t)\mbox{\rm d}W(t), \\[2mm]
P_{3}^{\varepsilon }(1) & = & 0.%
\end{array}%
\right.
\end{equation*}%
}
\end{remark}

\begin{remark}
\textsl{The necessary condition of near optimal controls are derived in
terms of the near maximum condition in an integral form. It is well known
that, for exact optimality, the integral form and the pointwise form of the
maximum condition are equivalent, however it is certainly not the case for
near optimality. }
\end{remark}

\begin{remark}
\textsl{In the work of Bahlali, Khelfallah, Mezerdi \cite{BKB}, they also
considered the near optimal control problem for general FBSDEs where the
diffusion term doesn't contain control variable. However, inspired by this
paper, we have noticed that, if }$\sigma $\textsl{\ contains control
variable, then this problem becomes more difficult. This topic will be
carried out as our future publication.}
\end{remark}

\section{Example}

\label{sect:4}

We now validate our theoretical results of Section \ref{sect:3} by looking
an example which is modified from Zhou \cite{Z1998}. Observe that the FBSDEs
considered in this paper are linear, it is possible to implement our
principles directly.

\begin{example}[\textbf{Necessary condition}]
\label{e1}\textsl{Let the admissible control domain $\Gamma =\left[ 0,1%
\right] .$ Consider the following $\varepsilon $-optimal control problem%
\begin{equation*}
\min\limits_{u(\cdot )\in \mathcal{U}_{ad}\left[ 0,1\right] }J(u(\cdot )),
\end{equation*}%
where
\begin{equation}
J(u(\cdot ))=\mathbb{E}\left[ \int_{0}^{1}u(t)\mbox{\rm d}t+\frac{\sqrt{2}}{2%
}x^{2}(1)+x(1)-y(0)\right] ,  \label{4.2}
\end{equation}%
with
\begin{equation}
\left\{
\begin{array}{rcl}
\mbox{\rm d}x(t) & = & u(t)\mbox{\rm d}W(t), \\
\mbox{\rm d}y(t) & = & -\left( 1+\sqrt{2}\right) u(t)\mbox{\rm d}t+z(t)%
\mbox{\rm d}W(t), \\
x(0) & = & 0,\qquad y(1)=x(1).%
\end{array}%
\right.  \label{4.3}
\end{equation}%
Set }$\theta _{0}^{\varepsilon }=\frac{\sqrt{2}}{2}.$ \textsl{For a given
admissible triple $(x^{\varepsilon }(\cdot ),y^{\varepsilon }(\cdot
),u^{\varepsilon }(\cdot )),$ the corresponding first and second adjoint
equations are%
\begin{equation}
\left\{
\begin{array}{rcl}
\mbox{\rm d}p^{\varepsilon }(t) & = & k^{\varepsilon }(t)\mbox{\rm d}W(t),
\\
\mbox{\rm d}q^{\varepsilon }(t) & = & 0, \\
p^{\varepsilon }(1) & = & x^{\varepsilon }(1),\qquad q^{\varepsilon }(0)=%
\frac{\sqrt{2}}{2},%
\end{array}%
\right.  \label{4.4}
\end{equation}%
and
\begin{equation}
\left\{
\begin{array}{rcl}
\mbox{\rm d}P_{1}^{\varepsilon }(t) & = & Q_{1}^{\varepsilon }(t)\mbox{\rm d}%
W(t), \\
P_{1}^{\varepsilon }(1) & = & 1,%
\end{array}%
\right.  \label{4.5}
\end{equation}%
respectively. Obviously, by the uniqueness of equations (\ref{4.4}) and (\ref%
{4.5}), we derive
\begin{equation*}
\left\{
\begin{array}{lll}
p^{\varepsilon }(t) & = & u^{\varepsilon }(t)W(t), \\
k^{\varepsilon }(t) & = & u^{\varepsilon }(t), \\
q^{\varepsilon }(t) & = & \frac{\sqrt{2}}{2},~t\in \left[ 0,1\right] ,%
\end{array}%
\right.
\end{equation*}%
and
\begin{equation*}
\left\{
\begin{array}{rcl}
P_{1}^{\varepsilon }(t) & = & 1, \\
Q_{1}^{\varepsilon }(t) & = & 0,~t\in \left[ 0,1\right] ,%
\end{array}%
\right.
\end{equation*}%
respectively. On the other hand, Theorem \ref{thm:4} gives
\begin{equation*}
\mathbb{E}\left\{ \int_{0}^{1}\left[ \frac{1}{2}(u(t))^{2}+u(t)(k^{%
\varepsilon }(t)-u^{\varepsilon }(t)-1)+\frac{1}{2}(u^{\varepsilon
}(t))^{2}-k^{\varepsilon }(t)u^{\varepsilon }(t)+u^{\varepsilon }(t)\right] %
\mbox{\rm d}t\right\} \geq -C\varepsilon ^{\beta }.
\end{equation*}%
Hence a simple calculation shows that if
\begin{equation}
u^{\varepsilon }(t)+1-k^{\varepsilon }(t)\in \Gamma ,  \label{4.7}
\end{equation}%
then, we get
\begin{equation}
\mathbb{E}\left[ \int_{0}^{1}(k^{\varepsilon }(t)-1)^{2}\mbox{\rm d}t\right]
\leq C\varepsilon ^{\beta }.  \label{4.8}
\end{equation}%
The above condition reveals the \textquotedblleft minimum\textquotedblright\
qualification for the pair $(x^{\varepsilon }(\cdot ),u^{\varepsilon }(\cdot
))$ to be $\varepsilon $-optimal. Actually, $u^{\varepsilon
}(t)=1-\varepsilon ^{\frac{1}{2}}$ is one of the candidates for $\varepsilon
$-optimal. Indeed, if we choose $u^{\varepsilon }(t)=1-\varepsilon ^{\frac{1%
}{2}}$ with the corresponding state
\begin{equation*}
x^{\varepsilon }(t)=(1-\varepsilon ^{\frac{1}{2}})W(t),~t\in \left[ 0,1%
\right] ,
\end{equation*}%
then the solutions of\ first order adjoint equations are
\begin{equation*}
(p^{\varepsilon }(t),k^{\varepsilon }(t))=\Big(\Big(1-\varepsilon ^{\frac{1}{%
2}}\Big)W(t),1-\varepsilon ^{\frac{1}{2}}\Big).
\end{equation*}%
Obviously, (\ref{4.7}) and (\ref{4.8}) are fulfilled. }
\end{example}

\section{Concluding Remarks}

\label{sect:5}

In this article, by Ekeland's principle, a spike variation, some dedicated
estimates and reduction method, we have established necessary condition for
near optimal controls to stochastic recursive optimization problems in terms
of a small parameter $\varepsilon >0$. In particular, we solve the problems
posed in \cite{HLW} (Huang, Li and Wang \textit{Near optimal control
problems for linear forward-backward stochastic systems}, Automatica 46
(2010), 397-404) Page 402 for control domain which is not necessarily convex
and diffusion term containing control variable. This result is partially
based on the work from \cite{BKB, E, HLW, W, Y, Z1998} etc. Our results
extends that of Zhou's \cite{Z1998} with second order adjoint equations in
the setup of FBSDEs. Hopefully, the theoretical result obtained in this
paper may inspire some real applications in finance and economics.

\appendix

\section*{Appendix}

\section{The Proof of Theorem \protect\ref{thm:4}}

\label{sect:6}

To establish the necessary condition, we need the following results mainly
from Lemma 2.1, Lemma 2.2, Lemma 3.1 and Lemma 3.2 in \cite{HLW} (note that $%
\left\vert \theta _{0}^{\varepsilon }\right\vert ^{2}+\mathbb{E}\left\vert
\theta _{1}^{\varepsilon }\right\vert ^{2}=1$, $1\geq \theta
_{0}^{\varepsilon }\geq 0$ for any fixed $\varepsilon >0$ which don't change
these results$)$. For simplicity, we omit the superscript $\varepsilon $.

\begin{lemma}
\label{lem:13} There exists a constant $C>0$ such that for any $\alpha \geq
0 $ and any $u( \cdot ) \in \mathcal{U}_{ad}\left[ 0,1\right] ,$%
\begin{equation*}
\begin{array}{c}
\mathbb{E}\left[ \sup\limits_{0\leq t\leq 1}\left\vert x( t) \right\vert
^{\alpha }\right] \leq C, \quad \mathbb{E}\left[ \sup\limits_{0\leq t\leq
1}\left\vert y( t) \right\vert ^{\alpha }\right] \leq C.%
\end{array}%
\end{equation*}
\end{lemma}

\begin{lemma}
\label{lem:14} There exists a constant $C>0$ such that
\begin{equation*}
\mathbb{E}\left[ \sup\limits_{0\leq t\leq 1}\left\vert q(t)\right\vert
^{2}+\sup\limits_{0\leq t\leq 1}\left\vert p(t)\right\vert
^{2}+\int_{0}^{1}\left\vert k(t)\right\vert ^{2}\mbox{\rm d}t\right] \leq C,
\end{equation*}%
where $C$ is independent of $(x(\cdot ),y(\cdot ),z(\cdot )).$
\end{lemma}

\begin{lemma}
There exists a constant $C>0$ such that
\begin{equation*}
\mathbb{E}\left[ \sup\limits_{0\leq t\leq 1}\left\vert P_{1}(t)\right\vert
^{2}+\int_{0}^{1}\left\vert Q_{1}(t)\right\vert ^{2}\mbox{\rm d}t\right]
\leq C,
\end{equation*}%
where $C$ is independent of $(x(\cdot ),y(\cdot ),z(\cdot )).$
\end{lemma}

\proof%
Applying Itô's formula to $\left\vert P_{3}(t)\right\vert ^{2},$ we have
\begin{equation*}
\left\vert P_{3}\left( t\right) \right\vert ^{2}+\mathbb{E}^{\mathcal{F}_{t}}%
\left[ \int_{t}^{1}\left\vert Q_{3}\left( s\right) \right\vert ^{2}%
\mbox{\rm
d}s\right] \leq C^{\prime }\mathbb{E}\int_{t}^{1}\left\vert P_{3}\left(
s\right) \right\vert ^{2}\mbox{\rm d}s+2\int_{t}^{1}\left\vert
l_{yy}\right\vert ^{2}\mbox{\rm d}s.
\end{equation*}%
By Burkholder-Davis-Gundy's inequality and Gronwall inequality, there is a
constant $C$ such that,
\begin{equation*}
\mathbb{E}\left[ \sup\limits_{0\leq t\leq 1}\left\vert P_{3}\left( t\right)
\right\vert ^{2}+\int_{t}^{1}\left\vert Q_{3}\left( s\right) \right\vert ^{2}%
\mbox{\rm d}s\right] \leq C.
\end{equation*}%
The same method to deal with $P_{2}\left( t\right) ,$ and $P_{1}\left(
t\right) ,$ we get the desired result.
\endproof%

\begin{lemma}
\label{lem:15} For any $\tau \geq 0$ and $0<\beta <1$ satisfying $\tau \beta
<1,$ there is a positive constant $C>0$ such that for any $u(\cdot )$ and $%
u^{\prime }(\cdot )\in \mathcal{U}_{ad}[0,1]$ along with the corresponding
trajectories $(x(\cdot ),y(\cdot ),z(\cdot ))$ and $(x^{\prime }(\cdot
),y^{\prime }(\cdot ),z^{\prime }(\cdot ))$, it follows that
\begin{equation*}
\left\{
\begin{array}{c}
\mathbb{E}\left[ \sup\limits_{0\leq t\leq 1}\left\vert x(t)-x^{\prime
}(t)\right\vert ^{2\tau }\right] \leq Cd(u(\cdot ),u^{\prime }(\cdot
))^{\tau \beta }, \\
\mathbb{E}\left[ \sup\limits_{0\leq t\leq 1}\left\vert y(t)-y^{\prime
}(t)\right\vert ^{2\tau }\right] \leq Cd(u(\cdot ),u^{\prime }(\cdot
))^{\tau \beta }.%
\end{array}%
\right.
\end{equation*}
\end{lemma}

\begin{lemma}
\label{lem:16} Assume (H1)-(H2) hold. For any $1<\tau <2$ and $0<\beta <1$
satisfying $(1+\beta )\tau <2,$ there is a constant $C$ such that for any $%
u(\cdot )$ and $u^{\prime }(\cdot )\in \mathcal{U}_{ad}[0,1]$ along with the
corresponding trajectories $\eta \left( \cdot \right) =(x(\cdot ),y(\cdot
),z(\cdot ))$, $\eta ^{\prime }\left( \cdot \right) =(x^{\prime }(\cdot
),y^{\prime }(\cdot ),z^{\prime }(\cdot ))$, and solutions $(p(\cdot
),q(\cdot ),k(\cdot )),(p^{\prime }(\cdot ),q^{\prime }(\cdot ),k^{\prime
}(\cdot ))$ of the corresponding adjoint equations, it holds that
\begin{equation*}
\left\{
\begin{array}{l}
\mathbb{E}\left[ \sup\limits_{0\leq t\leq 1}\left\vert q(t)-q^{\prime
}(t)\right\vert ^{\tau }\right] \leq Cd(u(\cdot ),u^{\prime }(\cdot ))^{%
\frac{\tau \beta }{2}}, \\[2mm]
\mathbb{E}\left[ \int_{0}^{1}\left( \left\vert p(t)-p^{\prime
}(t)\right\vert ^{\tau }+\left\vert k(t)-k^{\prime }(t)\right\vert ^{\tau
}\right) \mbox{\rm d}t\right] \leq Cd(u(\cdot ),u^{\prime }(\cdot ))^{\frac{%
\tau \beta }{2}}, \\[2mm]
\mathbb{E}\left[ \int_{0}^{1}\left( \left\vert P_{i}(t)-P_{i}^{\prime
}(t)\right\vert ^{\tau }+\left\vert Q_{i}(t)-Q_{i}^{\prime }(t)\right\vert
^{\tau }\right) \mbox{\rm d}t\right] \leq Cd(u(\cdot ),u^{\prime }(\cdot ))^{%
\frac{\tau \beta }{2}},~i=1,2,3.%
\end{array}%
\right.
\end{equation*}
\end{lemma}

\proof%
We are going to prove the third assertion. Note that $(\bar{P}_{3}(t),\bar{Q}%
_{3}(t))=(P_{3}(t)-P_{3}^{\prime }(t),Q_{3}(t)-Q_{3}^{\prime }(t))$
satisfies the following BSDEs%
\begin{equation*}
\left\{
\begin{array}{rcl}
-\text{d}\bar{P}_{3}(t) & = & \left[ -2b(t)\bar{P}_{3}(t)-\theta
_{0}(l_{yy}(t,x(t),y(t),u(t))-l_{yy}(t,x^{\prime }(t),y^{\prime
}(t)),u^{\prime }(t))\right] \mbox{\rm d}t \\[2mm]
&  & -Q_{3}(t)\mbox{\rm d}W(t), \\[2mm]
\bar{P}_{3}(1) & = & 0.%
\end{array}%
\right.
\end{equation*}%
Set $\rho _{3}(\cdot )$ to be the following linear SDEs:%
\begin{equation}
\left\{
\begin{array}{rcl}
\text{d}\rho _{3}(t) & = & \left[ 2b(t)\rho _{3}(t)+\left\vert \bar{P}%
_{3}(t)\right\vert ^{\tau -1}\text{sgn}(\bar{P}_{3}(t))\right] \mbox{\rm d}%
t+\left\vert \bar{Q}_{3}(t)\right\vert ^{\tau -1}\text{sgn}(\bar{Q}_{3}(t))%
\mbox{\rm d}W(t), \\
\rho _{3}(0) & = & 0,%
\end{array}%
\right.  \label{6.6}
\end{equation}%
It is easy to check that (\ref{6.6}) admit a unique solution, and the
following estimate can be obtained by Cauchy-Schwartz's inequality
\begin{equation}
\mathbb{E}\left[ \sup\limits_{0\leq t\leq 1}\left\vert \rho
_{3}(t)\right\vert ^{\gamma }\right] \leq C\mathbb{E}\left[
\int_{0}^{1}\left( \left[ \left\vert \bar{P}_{3}(t)\right\vert ^{\tau
}+\left\vert \bar{Q}_{3}(t)\right\vert ^{\tau }\right] \right) \mbox{\rm d}t%
\right] ,  \label{6.7}
\end{equation}%
where $\gamma >2$ and $\frac{1}{\gamma }+\frac{1}{\tau }=1$.

Applying Itô's formula to $\bar{P}_{3}(\cdot )\rho _{3}(\cdot )$ on $\left[
0,1\right] ,$ we have
\begin{align}
& \;\mathbb{E}\left[ \int_{0}^{1}\left( \left\vert \bar{P}_{3}(t)\right\vert
^{\tau }+\left\vert Q_{3}(t)\right\vert ^{\tau }\right) \mbox{\rm d}t\right]
\notag  \label{6.8} \\
=& \;\mathbb{E}\left[ \int_{0}^{1}\left( \rho _{3}(t)\theta
_{0}(l_{yy}(t,x(t),y(t),u(t))-l_{yy}(t,x^{\prime }(t),y^{\prime
}(t),u^{\prime }(t)))\right) \mbox{\rm d}t\right]  \notag \\
\leq & \;C\left( \mathbb{E}\int_{0}^{1}\left( \left\vert
(l_{yy}(t,x(t),y(t),u(t))-l_{yy}(t,x^{\prime }(t),y^{\prime }(t),u^{\prime
}(t)))\right\vert ^{\tau }\right) \mbox{\rm d}t\right) ^{\frac{1}{\tau }%
}\left( \mathbb{E}\int_{0}^{1}\left\vert \rho _{3}(t)\right\vert ^{\gamma }%
\mbox{\rm
d}t\right) ^{\frac{1}{\gamma }}.
\end{align}%
Substituting (\ref{6.7}) into (\ref{6.8}), we get
\begin{eqnarray*}
&&\mathbb{E}\left[ \int_{0}^{1}\left( \left\vert \bar{P}_{3}(t)\right\vert
^{\tau }+\left\vert Q_{3}(t)\right\vert ^{\tau }\right) \mbox{\rm d}t\right]
\\
&\leq &C\mathbb{E}\left[ \int_{0}^{1}\left( \theta _{0}\left\vert
(l_{yy}(t,x(t),y(t),u(t))-l_{yy}(t,x^{\prime }(t),y^{\prime }(t),u^{\prime
}(t)))\right\vert ^{\tau }\right) \mbox{\rm d}t\right] .
\end{eqnarray*}%
From (H1), it follows that
\begin{align*}
& \left[ \int_{0}^{1}\left( \theta _{0}\left\vert
l_{yy}(t,x(t),y(t),u(t))-l_{yy}(t,x^{\prime }(t),y^{\prime }(t),u^{\prime
}(t))\right\vert ^{\tau }\right) \mbox{\rm d}t\right] \\
\leq & \;\mathbb{E}\left[ \int_{0}^{1}\left( \left\vert
l_{yy}(t,x(t),y(t),u(t))-l_{yy}(t,x(t),y(t),u^{\prime }(t))\right\vert
^{\tau }\mathcal{X}_{u(t)=u^{^{\prime }}(t)}\right) \mbox{\rm d}t\right] \\
& +\mathbb{E}\left[ \int_{0}^{1}\left( \left\vert
l_{yy}(t,x(t),y(t),u^{\prime }(t))-l_{yy}(t,x^{\prime }(t),y^{\prime
}(t),u^{\prime }(t))\right\vert ^{\tau }\right) \mbox{\rm d}t\right] \\
\leq & \;C\mathbb{E}\left[ \left( \int_{0}^{1}\left( \left\vert
l_{yy}(t,x(t),y(t),u(t))-l_{yy}(t,x(t),y(t),u^{\prime }(t))\right\vert
^{2}\right) \mbox{\rm d}t\right) ^{\frac{\tau }{2}}\int_{0}^{1}d(u(t),u^{%
\prime }(t))^{1-\frac{\tau }{2}}\mbox{\rm d}t\right] \\
& +C\mathbb{E}\left[ \int_{0}^{1}\left( \left\vert x(t)-x^{\prime
}(t)\right\vert ^{\tau }+\left\vert y(t)-y^{\prime }(t)\right\vert ^{\tau
}\right) \mbox{\rm d}t\right] \\
\leq & \;Cd(u(t),u^{\prime }(t))^{\frac{\tau \beta }{2}}.
\end{align*}%
Combining (\ref{6.8}) with the above inequality, the result for $i=3$ holds
immediately.

We proceed to estimate the case, $i=2$. Similarly, we define the following
SDEs:
\begin{equation*}
\left\{
\begin{array}{rcl}
\text{d}\rho _{2}(t) & = & \left[ (b(t)-A(t))\rho _{2}(t)+\left\vert \bar{P}%
_{2}(t)\right\vert ^{\tau -1}\text{sgn}(\bar{P}_{2}(t))\right] \mbox{\rm d}t
\\
&  & +\left[ c(t)\rho _{2}(t)+\left\vert \bar{Q}_{2}(t)\right\vert ^{\tau -1}%
\text{sgn}(\bar{Q}_{2}(t))\right] \mbox{\rm d}W(t), \\
\rho _{2}(0) & = & 0.%
\end{array}%
\right.
\end{equation*}%
Applying Itô's formula to $\bar{P}_{2}(\cdot )\rho _{2}(\cdot )$ on $\left[
0,1\right] ,$ we have
\begin{align*}
& \;\mathbb{E}\left[ \int_{0}^{1}\left( \left\vert \bar{P}_{2}(t)\right\vert
^{\tau }+\left\vert Q_{2}(t)\right\vert ^{\tau }\right) \mbox{\rm d}t\right]
\\
=& \;\mathbb{E}\left[ \int_{0}^{1}\left( \rho _{2}(t)\left[ \theta
_{0}(l_{xy}(t,x(t),y(t),u(t))-l_{xy}(t,x^{\prime }(t),y^{\prime
}(t),u^{\prime }(t)))+a(t)\bar{P}_{3}(t)\right] \right) \mbox{\rm d}t\right]
.
\end{align*}%
By Cauchy-Schwartz's inequality, we obtain
\begin{align*}
& \;\mathbb{E}\left[ \int_{0}^{1}\left( \left\vert \bar{P}_{2}(t)\right\vert
^{\tau }+\left\vert Q_{2}(t)\right\vert ^{\tau }\right) \mbox{\rm d}t\right]
\\
\leq & \;C\mathbb{E}\left[ \int_{0}^{1}\left\vert a(t)\bar{P}%
_{3}(t)\right\vert ^{\tau }\mbox{\rm d}t\right] \\
& +C\mathbb{E}\left[ \int_{0}^{1}\left( \theta _{0}\left\vert
l_{xy}(t,x(t),y(t),u(t))-l_{xy}(t,x^{\prime }(t),y^{\prime }(t),u^{\prime
}(t))\right\vert ^{\tau }\right) \mbox{\rm d}t\right] \\
\leq & \;Cd(u(\cdot ),u^{\prime }(\cdot ))^{\frac{\tau \beta }{2}}.
\end{align*}%
Analogously, we define the following SDEs:
\begin{equation*}
\left\{
\begin{array}{lll}
\text{d}\rho _{1}(t) & = & \left[ (2A(t)+C(t)^{2})\rho _{1}(t)+\left\vert
\bar{P}_{1}(t)\right\vert ^{\tau -1}\text{sgn}(\bar{P}_{1}(t))\right] %
\mbox{\rm d}t \\
&  & +\left[ 2C(t)\rho _{1}(t)+\left\vert \bar{Q}_{1}(t)\right\vert ^{\tau
-1}\text{sgn}(\bar{Q}_{1}(t))\right] \mbox{\rm d}W(t), \\
\rho _{1}(0) & = & 0.%
\end{array}%
\right.
\end{equation*}%
Repeating the method used above, we have
\begin{equation*}
\mathbb{E}\left[ \int_{0}^{1}\left( \left\vert \bar{P}_{1}(t)\right\vert
^{\tau }+\left\vert Q_{1}(t)\right\vert ^{\tau }\right) \mbox{\rm d}t\right]
\leq Cd(u(t),u^{\prime }(t))^{\frac{\tau \beta }{2}}.
\end{equation*}%
The proof is complete.
\endproof%

The proof of Theorem \ref{thm:4} will be accomplished step by step. As the
reduction method developed by Yong \cite{Y} and Wu \cite{W}, independently,
we adopt the method by Yong \cite{Y} to derive the first and second adjoint
equations and the idea by Wu \cite{W} to deal with unbounded control problem
together.

\bigskip

\noindent \textit{Proof of Theorem \ref{thm:4}.}

\textbf{Step 1} \textbf{(The bounded control domains).}

When $(x(\cdot ),y(\cdot ))$ is regarded as the state process and $(z(\cdot
),u(\cdot ))$ as the control process, we consider the following initial
value problem for a control system of SDEs:%
\begin{equation}
\left\{
\begin{array}{rcl}
\mbox{\rm d}x(t) & = & \left[ A(t)x(t)+B(t)u(t)\right] \mbox{\rm d}t+\left[
C(t)x(t)+D(t)u(t)\right] \mbox{\rm d}W(t), \\[2mm]
-\mbox{\rm d}y(t) & = & \left[ a(t)x(t)+b(t)y(t)+c(t)u(t)\right] \text{d}%
t-z(t)\mbox{\rm d}W(t), \\[2mm]
x(0) & = & x_{0},\qquad y(0)=y_{0},%
\end{array}%
\right.  \label{6.16}
\end{equation}

\noindent Clearly, it is easy to check that, for any $(z(\cdot ),u(\cdot
))\in \mathcal{M}^{2}(0,1;\mathbb{R})\times \mathcal{U}_{ad}\left[ 0,1\right]
,$ $y_{0}\in \mathbb{R}$, there exists a unique strong solution
\begin{equation*}
(x(\cdot ),y(\cdot ))\equiv (x(\cdot ,z(\cdot ),u(\cdot )),y(\cdot ,z(\cdot
),u(\cdot )))\in \mathcal{S}^{2}(0,1;\mathbb{R})\times \mathcal{S}^{2}(0,1;%
\mathbb{R})
\end{equation*}%
to (\ref{6.16}) depending on $(z(\cdot ),u(\cdot ))$. Next, we regard the
original terminal condition as the terminal state constraint:
\begin{equation}
y(1)=Mx(1).  \label{6.17}
\end{equation}%
Since $\mathbb{R},$ $\mathcal{M}^{2}(0,1;\mathbb{R})$ are all unbounded,
Thus, we adopt a convergence technique developed by Wu \cite{W}.

Let $y_{0}$, $z\left( \cdot \right) $ take value in $\mathbb{M},$ $\mathbb{%
N\subset R},$ and $\mathbb{M}$ be convex. Moreover, $\mathbb{M},$ $\mathbb{N}
$ are all bounded. Let $\mathcal{A}$ be the set of all 3-triples $%
(y_{0},z(\cdot ),u(\cdot ))\in \mathbb{M}\times \mathcal{M}^{2}(0,1;\mathbb{N%
})\times \mathcal{U}_{ad}\left[ 0,1\right] $ such that the unique
corresponding state process $(x(\cdot ),y(\cdot ))$ satisfies the constraint
(\ref{6.17}). Note that, for any $u(\cdot )\in \mathcal{U}_{ad}[0,1]$, there
exists a unique $(y_{0},z(\cdot ))\in \mathbb{R}\times \mathcal{M}^{2}(0,1;%
\mathbb{R})$ such that state equation (\ref{2.1}) admits a unique state
process $(x(\cdot ),y(\cdot ))\in \mathcal{S}^{2}(0,1;\mathbb{R})\times
\mathcal{M}^{2}(0,1;\mathbb{R})$ satisfying the state constraint (\ref{6.17}%
). Hence, (H1) implies $\mathcal{A\neq \phi }.$ The cost functional is given
by
\begin{equation*}
J(y_{0},z(\cdot ),u(\cdot ))=\mathbb{E}\left[ \int_{0}^{1}l(t,x(t),y(t),u(t))%
\mbox{\rm d}t+\phi (x(1))+\gamma (y(0))\right] .
\end{equation*}%
We state the following problem.

\noindent \textbf{Problem $(\tilde{C}^{\varepsilon })$}. Find $%
(y_{0}^{\varepsilon },z^{\varepsilon }(\cdot ),u^{\varepsilon }(\cdot ))\in
\mathcal{A}$, such that
\begin{equation*}
J(y_{0}^{\varepsilon },z^{\varepsilon }(\cdot ),u^{\varepsilon }(\cdot
))=\inf\limits_{(y_{0},z(\cdot ),u(\cdot ))\in \mathcal{A}}J(y_{0},z(\cdot
),u(\cdot ))+\varepsilon .
\end{equation*}%
We, respectively, refer to $(y_{0}^{\varepsilon },z^{\varepsilon }(\cdot
),u^{\varepsilon }(\cdot ))$ as a near optimal control process, to $%
(x^{\varepsilon }(\cdot ),y^{\varepsilon }(\cdot ))$ as the corresponding
near optimal state process, and to $(y_{0}^{\varepsilon },z^{\varepsilon
}(\cdot ),u^{\varepsilon }(\cdot ))$ as a near optimal 3-tuple of Problem $(%
\tilde{C}^{\varepsilon })$.

Problems $(C^{\varepsilon })$ is embedded into $(\tilde{C}^{\varepsilon })$.
Suppose that $\left( \tilde{y}_{0}^{\varepsilon },\tilde{z}^{\varepsilon
}(\cdot ),\tilde{u}^{\varepsilon }(\cdot )\right) $ is the near optimal
control of Problem $(\tilde{C}^{\varepsilon })$, clearly, we know that $%
\tilde{u}^{\varepsilon }(\cdot )$ is the near optimal control of Problem $%
(C^{\varepsilon }).$ The advantage of Problem $(\tilde{C}^{\varepsilon })$
is that one does not need much regularity/integrability on $z(\cdot )$ since
it is treated as part of a control process; the disadvantage is that one has
to treat terminal constraint (\ref{6.17}).

\begin{lemma}[\textbf{Ekeland Principle \protect\cite{E}}]
\label{lem:17} Let $(S,d)$ be a complete metric space and $\rho
:S\rightarrow R\cup \left\{ +\infty \right\} $ be a lower semicontinuous
function, bounded from below. If for each $\varepsilon >0,$ there exists $%
u^{\varepsilon }\in S$ such that $\rho (u^{\varepsilon })\leq
\inf\limits_{u\in S}\rho (u)+\varepsilon .$ Then for any $\lambda >0,$ there
exists $u^{\lambda }\in S$ such that
\begin{equation*}
\left\{
\begin{array}{ll}
\mbox{(i)} & \rho (u^{\lambda })\leq \rho (u^{\varepsilon }), \\
\mbox{(ii)} & d(u^{\lambda },u^{\varepsilon })\leq \lambda , \\
\mbox{(iii)} & \rho (u^{\lambda })\leq \rho (u)+\frac{\varepsilon }{\lambda }%
d(u,u^{\lambda }),\quad \text{for all }u\in S.%
\end{array}%
\right.
\end{equation*}
\end{lemma}

\noindent For $u,v$ in $\mathcal{U}_{ad}\left[ 0,1\right] $ or in $\mathcal{M%
}^{2}(0,1;\mathbb{R}),$ we define
\begin{equation*}
d(u,v)=\text{d}t\otimes P\left\{ (t,\omega )\in \left[ 0,1\right] \times
\Omega :u(t,\omega )\neq v(t,\omega )\right\} ,
\end{equation*}%
where d$t\otimes P$ is the product measure of the Lebesgue measure d$t$ with
the probability measure $P$. It is well known that $(\mathcal{U}_{ad}\left[
0,1\right] ,d)$ is a complete metric space (see \cite{YZ}). Then $\mathbb{%
R\times }\mathcal{M}^{2}(0,1;\mathbb{R})\times \mathcal{U}_{ad}\left[ 0,1%
\right] $ is a complete metric space under the following metric: for any $%
(y_{0},z(\cdot ),u(\cdot )),(\tilde{y}_{0},\tilde{z}(\cdot ),\tilde{u}(\cdot
))\in \mathcal{A},$
\begin{equation*}
d_{\mathcal{A}}(\theta (\cdot ),\tilde{\theta}(\cdot ))=\left[ \left\vert
y_{0}-\tilde{y}_{0}\right\vert ^{2}+d\left( z(\cdot ),\tilde{z}(\cdot
)\right) ^{2}+d(u(\cdot ),\tilde{u}(\cdot ))^{2}\right] ^{\frac{1}{2}},
\end{equation*}%
where $\theta (\cdot )=(y_{0},z(\cdot ),u(\cdot ))$ and $\tilde{\theta}%
(\cdot )=(\tilde{y}_{0},\tilde{z}(\cdot ),\tilde{u}(\cdot )),$ respectively.

By assumption (H1), it is easy to see that $J(y_{0},z(\cdot ),u(\cdot ))$ is
lower semicontinuous on $\mathcal{A}$. By virtue of Ekeland principle (Lemma %
\ref{lem:17}) with $\lambda =\varepsilon ^{\frac{2}{3}}$ (fixed $\varepsilon
>0$) there is an admissible 3-triple $(\tilde{y}_{0}^{\varepsilon },\tilde{z}%
^{\varepsilon }(\cdot ),\tilde{u}^{\varepsilon }(\cdot ))\in \mathcal{A}$
such that
\begin{equation}
d_{\mathcal{A}}((y_{0}^{\varepsilon },z^{\varepsilon }(\cdot
),u^{\varepsilon }(\cdot )),(\tilde{y}_{0}^{\varepsilon },\tilde{z}%
^{\varepsilon }(\cdot ),\tilde{u}^{\varepsilon }(\cdot )))\leq \varepsilon ^{%
\frac{2}{3}}  \label{6.18}
\end{equation}%
and
\begin{equation*}
\tilde{J}^{\varepsilon }((\tilde{y}_{0}^{\varepsilon },\tilde{z}%
^{\varepsilon }(\cdot ),\tilde{u}^{\varepsilon }(\cdot )))\leq \tilde{J}%
^{\varepsilon }(v(\cdot )),\text{ for any }v(\cdot )\in \mathcal{A},
\end{equation*}%
where
\begin{equation}
\tilde{J}^{\varepsilon }(v(\cdot ))=J(v(\cdot ))+\varepsilon ^{\frac{1}{3}%
}d_{\mathcal{A}}(v(\cdot ),\tilde{\theta}^{\varepsilon }(\cdot )),
\label{6.22}
\end{equation}%
which means that $(\tilde{y}_{0},\tilde{z}^{\varepsilon }(\cdot ),\tilde{u}%
^{\varepsilon }(\cdot ))$ is an optimal triple for the system (\ref{6.16})
with a new cost functional $\tilde{J}^{\varepsilon }$.

Let $(\tilde{y}_{0}^{\varepsilon },\tilde{z}^{\varepsilon }(\cdot ),\tilde{u}%
^{\varepsilon }(\cdot ))$ be an optimal 3-triple of Problem $(\tilde{C}%
^{\varepsilon })$ with new functional (\ref{6.22}), with the corresponding
optimal state process $(\tilde{x}^{\varepsilon }(\cdot ),\tilde{y}%
^{\varepsilon }(\cdot )).$ For any $\delta >0,$ we define, for any $\forall
(y_{0},z(\cdot ),u(\cdot ))\in \mathbb{M}\times \mathcal{M}^{2}(0,1;\mathbb{N%
})\times \mathcal{U}_{ad}\left[ 0,1\right] $,
\begin{equation*}
J^{\delta ,\varepsilon }(y_{0},z(\cdot ),u(\cdot ))=\left\{ \left[ (\tilde{J}%
^{\varepsilon }(y_{0},z(\cdot ),u(\cdot ))-\tilde{J}^{\varepsilon }(\tilde{y}%
_{0}^{\varepsilon },\tilde{z}^{\varepsilon }(\cdot ),\tilde{u}^{\varepsilon
}(\cdot ))+\delta )^{+}\right] ^{2}+\mathbb{E}\left\vert
y(1)-Mx(1)\right\vert ^{2}\right\} ^{\frac{1}{2}},
\end{equation*}%
where $(x(\cdot ),y(\cdot ))$ is the unique solution of (\ref{6.16}). Also,
it is clear that
\begin{align*}
& J^{\delta ,\varepsilon }(y_{0},z(\cdot ),u(\cdot ))>0,~\forall
(y_{0},z(\cdot ),u(\cdot ))\in \mathbb{M}\times \mathcal{M}^{2}(0,1;\mathbb{N%
})\times \mathcal{U}_{ad}\left[ 0,1\right] , \\
& J^{\delta ,\varepsilon }(\tilde{y}_{0}^{\varepsilon },\tilde{z}%
^{\varepsilon }(\cdot ),\tilde{u}^{\varepsilon }(\cdot ))=\delta \leq
\inf\limits_{(y_{0},z(\cdot ),u(\cdot ))\in \mathbb{M}\times \mathcal{M}%
^{2}(0,1;\mathbb{N})\times \mathcal{U}_{ad}\left[ 0,1\right] }J^{\delta
,\varepsilon }(y_{0},z(\cdot ),u(\cdot ))+\delta .
\end{align*}%
Hence, by Lemma \ref{lem:17}, there exists a 3-triple $(y_{0}^{\delta
,\varepsilon },z^{\delta ,\varepsilon }(\cdot ),u^{\delta ,\varepsilon
}(\cdot ))\in \mathbb{M}\times \mathcal{M}^{2}(0,1;\mathbb{N})\times
\mathcal{U}_{ad}\left[ 0,1\right] $ such that
\begin{equation}
\left\{
\begin{array}{l}
\text{(1) }J^{\delta ,\varepsilon }(y_{0}^{\delta ,\varepsilon },z^{\delta
,\varepsilon }(\cdot ),u^{\delta ,\varepsilon }(\cdot ))\leq J^{\delta
,\varepsilon }(\tilde{y}_{0}^{\varepsilon },\tilde{z}^{\varepsilon }(\cdot ),%
\tilde{u}^{\varepsilon }(\cdot ))=\delta , \\[2mm]
\text{(2) }\left\vert y_{0}^{\delta ,\varepsilon }-\tilde{y}%
_{0}^{\varepsilon }\right\vert ^{2}+d\left( z^{\delta ,\varepsilon }(\cdot
)-z(\cdot )\right) ^{2}+d(u^{\delta ,\varepsilon }(\cdot ),\tilde{u}%
^{\varepsilon }(\cdot ))^{2}\leq \delta , \\[2mm]
\text{(3) }-\sqrt{\delta }\left[ \left\vert y_{0}^{\delta ,\varepsilon
}-y_{0}\right\vert ^{2}+d\left( z^{\delta ,\varepsilon }(\cdot )-z(\cdot
)\right) ^{2}+d(u^{\delta ,\varepsilon }(\cdot ),u(\cdot ))^{2}\right] ^{%
\frac{1}{2}} \\[2mm]
\qquad \leq J^{\delta ,\varepsilon }(y_{0},z(\cdot ),u(\cdot ))-J^{\delta
,\varepsilon }(y_{0}^{\delta ,\varepsilon },z^{\delta ,\varepsilon }(\cdot
),u^{\delta ,\varepsilon }(\cdot )), \\[2mm]
\hfill \forall (y_{0},z(\cdot ),u(\cdot ))\in \mathbb{M}\times \mathcal{M}%
^{2}(0,1;\mathbb{N})\times \mathcal{U}_{ad}\left[ 0,1\right] .%
\end{array}%
\right.  \label{6.26}
\end{equation}%
Hence, $(y_{0}^{\delta ,\varepsilon },z^{\delta ,\varepsilon }(\cdot
),u^{\delta ,\varepsilon }(\cdot ))$ is a global minimum point of the
following penalized cost functional%
\begin{equation}
J^{\delta ,\varepsilon }(y_{0},z(\cdot ),u(\cdot ))+\sqrt{\delta }\left[
\left\vert y_{0}-y_{0}^{\delta ,\varepsilon }\right\vert ^{2}+d\left(
z^{\delta ,\varepsilon }(\cdot )-z(\cdot )\right) ^{2}+d(u^{\delta
,\varepsilon }(\cdot ),u(\cdot ))^{2}\right] ^{\frac{1}{2}}.  \label{6.27}
\end{equation}%
In other words, fix $\varepsilon >0,$ if we pose a penalized optimal control
problem with the state constraint (\ref{6.17}) and the cost functional (\ref%
{6.27}) , then $(y_{0}^{\delta ,\varepsilon },z^{\delta ,\varepsilon }(\cdot
),u^{\delta ,\varepsilon }(\cdot ))$ is an optimal 3-triple of the problem.
Note that this problem does not have state constraints, and the optimal
3-triple $(y_{0}^{\delta ,\varepsilon },z^{\delta ,\varepsilon }(\cdot
),u^{\delta ,\varepsilon }(\cdot ))$ approaches $(\tilde{y}_{0}^{\varepsilon
},\tilde{z}^{\varepsilon }(\cdot ),\tilde{u}^{\varepsilon }(\cdot ))$ as $%
\delta \rightarrow 0$. Let us turn back to the new cost functional
\begin{eqnarray}
\mathcal{J}^{\delta ,\varepsilon }(y_{0},z(\cdot ),u(\cdot )) &=&J^{\delta
,\varepsilon }(y_{0},z(\cdot ),u(\cdot ))  \notag \\
&&+\sqrt{\delta }\left[ \left\vert y_{0}-y_{0}^{\delta ,\varepsilon
}\right\vert ^{2}+d\left( z^{\delta ,\varepsilon }(\cdot )-z(\cdot )\right)
^{2}+d(u^{\delta ,\varepsilon }(\cdot ),u(\cdot ))^{2}\right] ^{\frac{1}{2}}.
\label{6.28}
\end{eqnarray}

\noindent Denote
\begin{equation*}
\begin{array}{cccc}
X\triangleq \left(
\begin{array}{c}
x \\
y%
\end{array}%
\right) , & v(\cdot )\triangleq \left(
\begin{array}{c}
z \\
u%
\end{array}%
\right) , & X_{0}\triangleq \left(
\begin{array}{c}
x_{0} \\
y_{0}%
\end{array}%
\right) , & X(1)\triangleq \left(
\begin{array}{c}
x(1) \\
y(1)%
\end{array}%
\right) ,%
\end{array}%
\end{equation*}%
\begin{align*}
\mathcal{B}(t,X,v(\cdot ))\triangleq & \left(
\begin{array}{c}
A(t)x(t)+B(t)u(t) \\
-a(t)x(t)-b(t)y(t)-c(t)u(t)%
\end{array}%
\right) , \\
\Sigma (t,X,v(\cdot ))\triangleq & \left(
\begin{array}{c}
C(t)x(t)+D(t)u(t) \\
z%
\end{array}%
\right) , \\
\Xi (X(0),X(1))\triangleq & \;\phi (x(1))+\gamma (y(0)), \\
\Pi (X(0),X(1))\triangleq & \left(
\begin{array}{c}
0 \\
y(1)-Mx(1)%
\end{array}%
\right) ,
\end{align*}%
and
\begin{align*}
\mathcal{H}& \triangleq \;\mathbb{R}^{2}\times L_{\mathcal{F}%
_{1}}^{2}(\Omega ;\mathbb{R}^{2})\equiv \mathbb{R}^{2}\times \mathcal{X}%
_{2}^{2}, \\
\mathcal{H}_{0}& \triangleq \;\mathbb{R}\times L_{\mathcal{F}%
_{1}}^{2}(\Omega ;\mathbb{R})\equiv \mathbb{R}\times \mathcal{X}_{1}^{2}.
\end{align*}%
Consequently,
\begin{align*}
\tilde{J}^{\varepsilon }(y_{0},z(\cdot ),u(\cdot ))=& \;\tilde{J}%
^{\varepsilon }(y_{0},v(\cdot )), \\
J^{\delta ,\varepsilon }(y_{0},z(\cdot ),u(\cdot ))=& \;J^{\delta
,\varepsilon }(y_{0},v(\cdot )).
\end{align*}%
Note that $\mathcal{H}$ and $\mathcal{H}_{0}$ are Hilbert spaces. We
identify $\mathcal{H}^{\ast }=\mathcal{H}$ and $\mathcal{H}_{0}^{\ast }=%
\mathcal{H}_{0}.$ Also
\begin{equation*}
\Xi :\mathcal{H}\rightarrow \mathbb{R},\qquad \Pi :\mathcal{H}\rightarrow
\mathcal{H}_{0}.
\end{equation*}%
The gradient of $D\Xi $ and the Hessian $D^{2}\Xi $ of $\Xi $ are defined as
follows:
\begin{align*}
D\Xi (X(0),X(1))=& (D_{X_{0}}\Xi (X(0),X(1)),D_{X_{1}}\Xi (X(0),X(1)))\in
\mathcal{L}(\mathcal{H};\mathbb{R})\equiv \mathcal{H}^{\ast }=\mathcal{H}, \\
D^{2}\Xi (X(0),X(1))=& \left(
\begin{array}{cc}
D_{X_{0}X_{0}}\Xi (X(0),X(1)) & D_{X_{0}X_{1}}\Xi (X(0),X(1)) \\
D_{X_{1}X_{0}}\Xi (X(0),X(1)) & D_{X_{1}X_{1}}\Xi (X(0),X(1))%
\end{array}%
\right) \in \mathcal{L}_{s}(\mathcal{H};\mathcal{H}),
\end{align*}%
where $\mathcal{L}(\mathcal{H}_{1};\mathcal{H}_{2})$ is the set of all
linear bounded operator from $\mathcal{H}_{1}$ to $\mathcal{H}_{2}$, and $%
\mathcal{L}_{s}(\mathcal{H};\mathcal{H})$ is the set of all linear bounded
self-adjoint operators from $\mathcal{H}$ to itself. We have%
\begin{eqnarray*}
\Xi _{X_{0}}(X(0),X(1)) &=&(0,\gamma _{y}(y_{0}))^{\top }\in \mathbb{R}^{2},
\\
\Xi _{X_{1}}(X(0),X(1)) &=&(\phi _{x}(x(1)),0)^{\top }\in \mathcal{X}%
_{1}^{2}, \\
\Xi _{X_{0}X_{0}}(X(0),X(1)) &=&\left(
\begin{array}{cc}
0 & 0 \\
0 & \gamma _{yy}(y_{0})%
\end{array}%
\right) \in \mathcal{S}^{2}, \\
\Xi _{X_{0}X_{1}}(X(0),X(1)) &=&\left(
\begin{array}{cc}
0 & 0 \\
0 & 0%
\end{array}%
\right) \in \mathcal{L}(\mathcal{X}_{1}^{2};\mathbb{R}^{2}), \\
\Xi _{X_{1}X_{0}}(X(0),X(1)) &=&\left(
\begin{array}{cc}
0 & 0 \\
0 & 0%
\end{array}%
\right) \in \mathcal{L}(\mathbb{R}^{2};\mathcal{X}_{1}^{2}), \\
\Xi _{X_{1}X_{1}}(X(0),X(1)) &=&\left(
\begin{array}{cc}
\phi _{xx}(x(1)) & 0 \\
0 & 0%
\end{array}%
\right) \in \mathcal{L}(\mathcal{X}_{2}^{2};\mathcal{X}_{2}^{2}).
\end{eqnarray*}
For
\begin{align*}
D\Pi (X(0),X(1))=& (D_{X_{0}}\Pi (X(0),X(1)),D_{X_{1}}\Pi (X(0),X(1)))\in
\mathcal{L}(\mathcal{H};\mathcal{H}_{0}), \\
D^{2}\Pi (X(0),X(1))=& (%
\begin{array}{cc}
D_{X_{0}X_{0}}\Pi (X(0),X(1)) & D_{X_{0}X_{1}}\Pi (X(0),X(1)) \\
D_{X_{1}X_{0}}\Pi (X(0),X(1)) & D_{X_{1}X_{1}}\Pi (X(0),X(1))%
\end{array}%
)\in \mathcal{L}(\mathcal{H};\mathcal{L}(\mathcal{H};\mathcal{H}_{0})).
\end{align*}%
Take any $\hat{\Phi}=(\hat{\Phi}_{0},\hat{\Phi}_{1})\in \mathcal{H}_{0}.$
Then,
\begin{equation*}
\left\langle \Pi (X(0),X(1)),\hat{\Phi}\right\rangle =\left\langle
y(1)-M(1)x(1),\hat{\Phi}\right\rangle .
\end{equation*}%
Thus,
\begin{align*}
D\Pi (X(0),X(1))\hat{\Phi}=& \;D\left[ \left\langle \Pi (X(0),X(1)),\hat{\Phi%
}\right\rangle \right] \\
=& \left( \left\langle \Pi (X(0),X(1)),\hat{\Phi}\right\rangle
_{X_{0}},\left\langle \Pi (X(0),X(1)),\hat{\Phi}\right\rangle _{X_{1}}\right)
\\
=& (\Pi _{X_{0}}(X(0),X(1))\hat{\Phi},\Pi _{X_{1}}(X(0),X(1))\hat{\Phi}),
\end{align*}%
with
\begin{align*}
\Pi _{X_{0}}(X(0),X(1))\hat{\Phi}=& (0,0), \\
\Pi _{X_{1}}(X(0),X(1))\hat{\Phi}=& (-M\hat{\Phi}_{1},\hat{\Phi}_{1}),
\end{align*}%
\begin{align*}
D^{2}\Pi (X(0),X(1))\hat{\Phi}=& \;D^{2}\left[ \left\langle \Pi (X(0),X(1)),%
\hat{\Phi}\right\rangle \right] \\
=& \left(
\begin{array}{cc}
D_{X_{0}X_{0}}\Pi (X(0),X(1))\hat{\Phi} & D_{X_{0}X_{1}}\Pi (X(0),X(1))\hat{%
\Phi} \\
D_{X_{1}X_{0}}\Pi (X(0),X(1))\hat{\Phi} & D_{X_{1}X_{1}}\Pi (X(0),X(1))\hat{%
\Phi}%
\end{array}%
\right) \in \mathcal{L}(\mathcal{H};\mathcal{H}),
\end{align*}%
and
\begin{align*}
D_{X_{0}X_{0}}\Pi (X(0),X(1))\hat{\Phi}=& \left(
\begin{array}{cc}
0 & 0 \\
0 & 0%
\end{array}%
\right) , \\
D_{X_{1}X_{0}}\Pi (X(0),X(1))\hat{\Phi}=& \left(
\begin{array}{cc}
0 & 0 \\
0 & 0%
\end{array}%
\right) , \\
D_{X_{0}X_{1}}\Pi (X(0),X(1))\hat{\Phi}=& \left(
\begin{array}{cc}
0 & 0 \\
0 & 0%
\end{array}%
\right) , \\
D_{X_{1}X_{1}}\Pi (X(0),X(1))\hat{\Phi}=& \left(
\begin{array}{cc}
0 & 0 \\
0 & 0%
\end{array}%
\right) .
\end{align*}%
We now construct spike variation. For$,$ $u(\cdot )\in \mathcal{U}_{ad}\left[
0,1\right] $, and any $0<\alpha <1,$ let $y_{0}\in \mathbb{M}$ \ such that $%
y_{0}^{\delta ,\varepsilon }+y_{0}\in \mathbb{M}$ . Define
\begin{equation*}
y_{0}^{\delta ,\varepsilon ,\alpha }=\alpha y_{0}+y_{0}^{\delta ,\varepsilon
},\text{\qquad }(z^{\delta ,\varepsilon ,\alpha }(t),u^{\delta ,\varepsilon
,\alpha }(t))=\left\{
\begin{array}{ll}
(z,,u), & t\in \left[ \tau ,\tau +\alpha \right] , \\
(z^{\delta ,\varepsilon }(t),u^{\delta ,\varepsilon }(t)), & \text{otherwise,%
}%
\end{array}%
\right.
\end{equation*}%
where $z\in \mathbb{N}$, $u\in \mathbb{U}$ are $\mathcal{F}_{\tau }$%
-measurable random variables, such that $\sup\limits_{\omega \in \Omega
}\left\vert u\left( \omega \right) \right\vert <+\infty $ and $%
\sup\limits_{\omega \in \Omega }\left\vert z\left( \omega \right)
\right\vert <+\infty $. Note that $y_{0}$ is a control independent of time
variable, so convex perturbation can be applied here.

Let $X^{\delta ,\varepsilon ,\alpha }(\cdot )$ be the state process (\ref%
{6.16}) corresponding to $(X_{0}^{\delta ,\varepsilon ,\alpha }\cdot
,v^{\delta ,\varepsilon ,\alpha }(\cdot ))$. Let $X_{1}^{\delta ,\varepsilon
,\alpha }(\cdot )$ and $X_{2}^{\delta ,\varepsilon ,\alpha }(\cdot )$ be,
respectively, the solutions to the following SDEs:
\begin{equation*}
\left\{
\begin{array}{lll}
\mbox{\rm d}X_{1}^{\delta ,\varepsilon ,\alpha }(t) & = & \mathcal{B}%
_{X}^{\delta ,\varepsilon }(t,\cdot )X_{1}^{\delta ,\varepsilon ,\alpha }(t)%
\mbox{\rm d}t+\left[ \Sigma _{X}^{\delta ,\varepsilon }(t,\cdot
)X_{1}^{\delta ,\varepsilon ,\alpha }(t)+\bigtriangleup \Sigma ^{\delta
,\varepsilon }(t,\cdot )I_{S_{\alpha }}\right] \mbox{\rm d}W(t), \\
X_{1}^{\delta ,\varepsilon ,\alpha }(t) & = & \sqrt{\alpha }y_{0},%
\end{array}%
\right.
\end{equation*}%
and
\begin{equation*}
\left\{
\begin{array}{lll}
\mbox{\rm d}X_{2}^{\delta ,\varepsilon ,\alpha }(t) & = & \mathcal{B}%
_{X}^{\delta ,\varepsilon }(t,\cdot )X_{2}^{\delta ,\varepsilon ,\alpha
}(t)+\bigtriangleup \mathcal{B}^{\delta ,\varepsilon }(t,\cdot )I_{S_{\alpha
}}+\frac{1}{2}\mathcal{B}_{XX}^{\delta ,\varepsilon }(t,\cdot )X_{1}^{\delta
,\varepsilon ,\alpha }(t)^{2}\mbox{\rm d}t \\[2mm]
&  & +\left[ \Sigma _{X}^{\delta ,\varepsilon }(t,\cdot )X_{2}^{\delta
,\varepsilon ,\alpha }(t)+\bigtriangleup \Sigma _{X}^{\delta ,\varepsilon
}(t,\cdot )X_{1}^{\delta ,\varepsilon ,\alpha }(t)I_{S_{\alpha }}+\frac{1}{2}%
\Sigma _{XX}^{\delta ,\varepsilon }(t,\cdot )X_{1}^{\delta ,\varepsilon
,\alpha }(t)^{2}\right] \mbox{\rm d}W(t), \\[2mm]
X_{2}^{\delta ,\varepsilon ,\alpha }(0) & = & 0,%
\end{array}%
\right.
\end{equation*}%
where $I_{S_{\alpha }}$ denotes the indicator function of the set $S_{\alpha
}$ and for any $X\in \mathbb{R}^{2}$.

Set
\begin{align*}
\mathcal{B}_{X}^{\delta ,\varepsilon }(t,\cdot )=& \;\mathcal{B}%
_{X}(t,X^{\delta ,\varepsilon }(t),v^{\delta ,\varepsilon }(t)), \\
\Sigma _{X}^{\delta ,\varepsilon }(t,\cdot )=& \;\Sigma _{X}(t,X^{\delta
,\varepsilon }(t),v^{\delta ,\varepsilon }(t)), \\
\bigtriangleup \mathcal{B}^{\delta ,\varepsilon }(t,\cdot )=& \;\mathcal{B}%
(t,X^{\delta ,\varepsilon }(t),v(t))-\mathcal{B}(t,X^{\delta ,\varepsilon
}(t),v^{\delta ,\varepsilon }(t)), \\
\bigtriangleup \Sigma ^{\delta ,\varepsilon }(t,\cdot )=& \;\Sigma
(t,X^{\delta ,\varepsilon }(t),v(t))-\Sigma (t,X^{\delta ,\varepsilon
}(t),v^{\delta ,\varepsilon }(t)), \\
\bigtriangleup \Sigma _{X}^{\delta ,\varepsilon }(t,\cdot )=& \;\Sigma
_{X}(t,X^{\delta ,\varepsilon }(t),v(t))-\Sigma _{X}(t,X^{\delta
,\varepsilon }(t),v^{\delta ,\varepsilon }(t)),
\end{align*}%
\begin{align*}
\mathcal{B}_{XX}^{\delta ,\varepsilon }(t,\cdot )X^{2}=& \;\left(
\begin{array}{c}
\left\langle \mathcal{B}_{XX}^{1,\delta ,\varepsilon }(t,\cdot
)X,X\right\rangle \\[2mm]
\left\langle \mathcal{B}_{XX}^{2,\delta ,\varepsilon }(t,\cdot
)X,X\right\rangle%
\end{array}%
\right) ,~i=1,2, \\
\mathcal{B}_{XX}^{i,\delta ,\varepsilon }(t,\cdot )=& \;\mathcal{B}%
_{XX}^{i}(t,X^{\delta ,\varepsilon }(t),v^{\delta ,\varepsilon }(t)),~i=1,2,
\\
\Sigma _{XX}^{\delta ,\varepsilon }(t,\cdot )=& \left(
\begin{array}{c}
\left\langle \Sigma _{XX}^{1,\delta ,\varepsilon }(t,\cdot )X,X\right\rangle
\\[2mm]
\left\langle \Sigma _{XX}^{2,\delta ,\varepsilon }(t,\cdot )X,X\right\rangle%
\end{array}%
\right) ,~i=1,2, \\
\Sigma _{XX}^{i,\delta ,\varepsilon }(t,\cdot )=& \;\Sigma
_{XX}^{i}(t,X^{\delta ,\varepsilon }(t),v^{\delta ,\varepsilon }(t)),~i=1,2.
\end{align*}%
The following results can be seen in Wu \cite{W}:
\begin{align*}
& \sup\limits_{0\leq t\leq 1}\mathbb{E}\left\vert X_{1}^{\delta ,\varepsilon
,\alpha }(t)\right\vert ^{2k}+\sup\limits_{0\leq t\leq 1}\mathbb{E}%
\left\vert X^{\delta ,\varepsilon ,\alpha }(t)-X^{\delta ,\varepsilon
}(t)\right\vert ^{2k}\leq C\alpha ^{k}, \\
& \sup\limits_{0\leq t\leq 1}\mathbb{E}\left\vert X_{2}^{\delta ,\varepsilon
,\alpha }(t)\right\vert ^{2k}+\sup\limits_{0\leq t\leq 1}\mathbb{E}%
\left\vert X^{\delta ,\varepsilon ,\alpha }(t)-X^{\delta ,\varepsilon
}(t)-X_{1}^{\delta ,\varepsilon ,\alpha }(t)\right\vert ^{2k}\leq C\alpha
^{2k}, \\
& \sup\limits_{0\leq t\leq 1}\mathbb{E}\left\vert X^{\delta ,\varepsilon
,\alpha }(t)-X^{\delta ,\varepsilon }(t)-X_{1}^{\delta ,\varepsilon ,\alpha
}(t)-X_{2}^{\delta ,\varepsilon ,\alpha }(t)\right\vert ^{2k}=o(\alpha
^{2k}).
\end{align*}%
Now from the last relation in (\ref{6.26}), we derive%
\begin{eqnarray*}
-\sqrt{\delta }\alpha \sqrt{2+\left\vert y_{0}\right\vert ^{2}} &\leq
&\;J^{\delta ,\varepsilon }(y_{0}^{\delta ,\varepsilon ,\alpha },z^{\delta
,\varepsilon ,\alpha }(\cdot ),u^{\delta ,\varepsilon ,\alpha }(\cdot
))-J^{\delta ,\varepsilon }(y_{0}^{\delta ,\varepsilon },z^{\delta
,\varepsilon }(\cdot ),u^{\delta ,\varepsilon }(\cdot )) \\
&=&\frac{J^{\delta ,\varepsilon }(y_{0}^{\delta ,\varepsilon ,\alpha
},z^{\delta ,\varepsilon ,\alpha }(\cdot ),u^{\delta ,\varepsilon ,\alpha
}(\cdot ))^{2}}{J^{\delta ,\varepsilon }(y_{0}^{\delta ,\varepsilon ,\alpha
},z^{\delta ,\varepsilon ,\alpha }(\cdot ),u^{\delta ,\varepsilon ,\alpha
}(\cdot ))-J^{\delta ,\varepsilon }(y_{0}^{\delta ,\varepsilon },z^{\delta
,\varepsilon }(\cdot ),u^{\delta ,\varepsilon }(\cdot ))} \\
&&-\frac{J^{\delta ,\varepsilon }(y_{0}^{\delta ,\varepsilon },z^{\delta
,\varepsilon }(\cdot ),u^{\delta ,\varepsilon }(\cdot ))^{2}}{J^{\delta
,\varepsilon }(y_{0}^{\delta ,\varepsilon ,\alpha },z^{\delta ,\varepsilon
,\alpha }(\cdot ),u^{\delta ,\varepsilon ,\alpha }(\cdot ))-J^{\delta
,\varepsilon }(y_{0}^{\delta ,\varepsilon },z^{\delta ,\varepsilon }(\cdot
),u^{\delta ,\varepsilon }(\cdot ))} \\
&=&\frac{\left[ (\tilde{J}^{\varepsilon }(y_{0}^{\delta ,\varepsilon ,\alpha
},z^{\delta ,\varepsilon ,\alpha }(\cdot ),u^{\delta ,\varepsilon ,\alpha
}(\cdot ))-\tilde{J}^{\varepsilon }(\tilde{y}_{0}^{\varepsilon },\tilde{z}%
^{\varepsilon }(\cdot ),\tilde{u}^{\varepsilon }(\cdot ))+\delta )^{+}\right]
^{2}}{J^{\delta ,\varepsilon }(y_{0}^{\delta ,\varepsilon ,\alpha
},z^{\delta ,\varepsilon ,\alpha }(\cdot ),u^{\delta ,\varepsilon ,\alpha
}(\cdot ))-J^{\delta ,\varepsilon }(y_{0}^{\delta ,\varepsilon },z^{\delta
,\varepsilon }(\cdot ),u^{\delta ,\varepsilon }(\cdot ))} \\
&&-\frac{\left[ (\tilde{J}^{\varepsilon }(y_{0}^{\delta ,\varepsilon
},z^{\delta ,\varepsilon }(\cdot ),u^{\delta ,\varepsilon }(\cdot ))-\tilde{J%
}^{\varepsilon }(\tilde{y}_{0}^{\varepsilon },\tilde{z}^{\varepsilon }(\cdot
),\tilde{u}^{\varepsilon }(\cdot ))+\delta )^{+}\right] ^{2}}{J^{\delta
,\varepsilon }(y_{0}^{\delta ,\varepsilon ,\alpha },z^{\delta ,\varepsilon
,\alpha }(\cdot ),u^{\delta ,\varepsilon ,\alpha }(\cdot ))-J^{\delta
,\varepsilon }(y_{0}^{\delta ,\varepsilon },z^{\delta ,\varepsilon }(\cdot
),u^{\delta ,\varepsilon }(\cdot ))} \\
&&+\frac{\mathbb{E}\left[ \Pi (X^{\delta ,\varepsilon ,\alpha }(0),X^{\delta
,\varepsilon ,\alpha }(1))^{2}-\Pi (X^{\delta ,\varepsilon }(0),X^{\delta
,\varepsilon }(1))^{2}\right] }{J^{\delta ,\varepsilon }(y_{0}^{\delta
,\varepsilon ,\alpha },z^{\delta ,\varepsilon ,\alpha }(\cdot ),u^{\delta
,\varepsilon ,\alpha }(\cdot ))-J^{\delta ,\varepsilon }(y_{0}^{\delta
,\varepsilon },z^{\delta ,\varepsilon }(\cdot ),u^{\delta ,\varepsilon
}(\cdot ))} \\
&=&\;\theta _{0}^{\delta ,\varepsilon ,\alpha }\left[ \tilde{J}^{\varepsilon
}(y_{0}^{\delta ,\varepsilon ,\alpha },z^{\delta ,\varepsilon ,\alpha
}(\cdot ),u^{\delta ,\varepsilon ,\alpha }(\cdot ))-\tilde{J}^{\varepsilon
}(y_{0}^{\delta ,\varepsilon },z^{\delta ,\varepsilon }(\cdot ),u^{\delta
,\varepsilon }(\cdot ))\right] \\
&&+\mathbb{E}\left\langle \binom{0}{\theta _{1}^{\delta ,\varepsilon ,\alpha
}},\Pi (X^{\delta ,\varepsilon ,\alpha }(0),X^{\delta ,\varepsilon ,\alpha
}(1))-\Pi (X^{\delta ,\varepsilon }(0),X^{\delta ,\varepsilon
}(1))\right\rangle \\
&=&(\theta _{0}^{\delta ,\varepsilon }+o(1))\left[ \tilde{J}^{\varepsilon
}(y_{0}^{\delta ,\varepsilon ,\alpha },z^{\delta ,\varepsilon ,\alpha
}(\cdot ),u^{\delta ,\varepsilon ,\alpha }(\cdot ))-\tilde{J}^{\varepsilon
}(y_{0}^{\delta ,\varepsilon },z^{\delta ,\varepsilon }(\cdot ),u^{\delta
,\varepsilon }(\cdot ))\right] \\
&&+\mathbb{E}\left\langle \left(
\begin{array}{c}
0 \\
\theta _{1}^{\delta ,\varepsilon }+o(1)%
\end{array}%
\right) ,\Pi (X^{\delta ,\varepsilon ,\alpha }(0),X^{\delta ,\varepsilon
,\alpha }(1))-\Pi (X^{\delta ,\varepsilon }(0),X^{\delta ,\varepsilon
}(1))\right\rangle ,
\end{eqnarray*}%
where%
\begin{align*}
\theta _{0}^{\delta ,\varepsilon ,\alpha }& =\frac{2}{J^{\delta ,\varepsilon
}(y_{0}^{\delta ,\varepsilon ,\alpha },z^{\delta ,\varepsilon ,\alpha
}(\cdot ),u^{\delta ,\varepsilon ,\alpha }(\cdot ))+J^{\delta ,\varepsilon
}(y_{0}^{\delta ,\varepsilon },z^{\delta ,\varepsilon }(\cdot ),u^{\delta
,\varepsilon }(\cdot ))} \\
& \times \Bigg \{\int_{0}^{1}\bigg[\beta (\tilde{J}(y_{0}^{\delta
,\varepsilon ,\alpha },z^{\delta ,\varepsilon ,\alpha }(\cdot ),u^{\delta
,\varepsilon ,\alpha }(\cdot ))-\tilde{J}(\tilde{y}_{0}^{\varepsilon },%
\tilde{z}^{\varepsilon }(\cdot ),\tilde{u}^{\varepsilon }(\cdot ))) \\
& +(1-\beta )(\tilde{J}(y_{0}^{\delta ,\varepsilon },z^{\delta ,\varepsilon
}(\cdot ),u^{\delta ,\varepsilon }(\cdot ))-\tilde{J}(\tilde{y}%
_{0}^{\varepsilon },\tilde{z}^{\varepsilon }(\cdot ),\tilde{u}^{\varepsilon
}(\cdot )))\bigg]\mbox{\rm d}\beta +\delta \Bigg \}^{+}, \\
\theta _{1}^{\delta ,\varepsilon ,\alpha }& =\frac{y^{\delta ,\varepsilon
,\alpha }(1)-Mx^{\delta ,\varepsilon ,\alpha }(1)+y^{\delta ,\varepsilon
}(1)-Mx^{\delta ,\varepsilon }(1)}{J^{\delta ,\varepsilon }(y_{0}^{\delta
,\varepsilon ,\alpha },z^{\delta ,\varepsilon ,\alpha }(\cdot ),u^{\delta
,\varepsilon ,\alpha }(\cdot ))+J^{\delta ,\varepsilon }(y_{0}^{\delta
,\varepsilon },z^{\delta ,\varepsilon }(\cdot ),u^{\delta ,\varepsilon
}(\cdot ))}, \\
\theta _{0}^{\delta ,\varepsilon }& =\frac{(\tilde{J}^{\varepsilon
}(y_{0}^{\delta ,\varepsilon },z^{\delta ,\varepsilon }(\cdot ),u^{\delta
,\varepsilon }(\cdot ))-\tilde{J}^{\varepsilon }(\tilde{y}_{0}^{\varepsilon
},\tilde{z}^{\varepsilon }(\cdot ),\tilde{u}^{\varepsilon }(\cdot ))+\delta
)^{+}}{J^{\delta ,\varepsilon }(y^{\delta ,\varepsilon }(0),z^{\delta
,\varepsilon }(\cdot ),u^{\delta ,\varepsilon }(\cdot ))}\in \left[ 0,1%
\right] , \\
\theta _{1}^{\delta ,\varepsilon }& =\frac{y^{\delta ,\varepsilon
}(1)-Mx^{\delta ,\varepsilon }(1)}{J^{\delta ,\varepsilon }(y^{\delta
,\varepsilon }(0),z^{\delta ,\varepsilon }(\cdot ),u^{\delta ,\varepsilon
}(\cdot ))}\in L_{\mathcal{F}_{1}}^{2}(\Omega ;\mathbb{R}).
\end{align*}%
On the other hand,
\begin{align*}
& \;\theta _{0}^{\delta ,\varepsilon ,\alpha }\left[ \tilde{J}^{\varepsilon
}(y_{0}^{\delta ,\varepsilon ,\alpha },z^{\delta ,\varepsilon ,\alpha
}(\cdot ),u^{\delta ,\varepsilon ,\alpha }(\cdot ))-\tilde{J}^{\varepsilon
}(y_{0}^{\delta ,\varepsilon },z^{\delta ,\varepsilon }(\cdot ),u^{\delta
,\varepsilon }(\cdot ))\right] \\
=& \;\theta _{0}^{\delta ,\varepsilon ,\alpha }\left[ J(y_{0}^{\delta
,\varepsilon ,\alpha },z^{\delta ,\varepsilon ,\alpha }(\cdot ),u^{\delta
,\varepsilon ,\alpha }(\cdot ))+\varepsilon ^{\frac{1}{3}}d_{\mathcal{A}%
}(y_{0}^{\delta ,\varepsilon ,\alpha },z^{\delta ,\varepsilon }(\cdot
),u^{\delta ,\varepsilon ,\alpha }(\cdot ),\tilde{\theta}^{\varepsilon
}(\cdot ))\right. \\
& \left. -J(y_{0}^{\delta ,\varepsilon },z^{\delta ,\varepsilon }(\cdot
),u^{\delta ,\varepsilon }(\cdot ))-\varepsilon ^{\frac{1}{3}}d_{\mathcal{A}%
}(y_{0}^{\delta ,\varepsilon },z^{\delta ,\varepsilon }(\cdot ),u^{\delta
,\varepsilon }(\cdot ),\tilde{\theta}^{\varepsilon }(\cdot ))\right] \\
=& \;\theta _{0}^{\delta ,\varepsilon ,\alpha }\left[ J(y_{0}^{\delta
,\varepsilon ,\alpha },z^{\delta ,\varepsilon ,\alpha }(\cdot ),u^{\delta
,\varepsilon ,\alpha }(\cdot ))-J(y_{0}^{\delta ,\varepsilon },z^{\delta
,\varepsilon }(\cdot ),u^{\delta ,\varepsilon }(\cdot ))\right] \\
& +\varepsilon ^{\frac{1}{3}}\theta _{0}^{\delta ,\varepsilon ,\alpha }\left[
d_{\mathcal{A}}(y_{0}^{\delta ,\varepsilon ,\alpha },z^{\delta ,\varepsilon
,\alpha }(\cdot ),u^{\delta ,\varepsilon ,\alpha }(\cdot ),\tilde{\theta}%
^{\varepsilon }(\cdot ))-d_{\mathcal{A}}(y_{0}^{\delta ,\varepsilon
},z^{\delta ,\varepsilon }(\cdot ),u^{\delta ,\varepsilon }(\cdot ),\tilde{%
\theta}^{\varepsilon }(\cdot ))\right] \\
\leq & \;\theta _{0}^{\delta ,\varepsilon ,\alpha }\left[ J(y_{0}^{\delta
,\varepsilon ,\alpha },z^{\delta ,\varepsilon ,\alpha }(\cdot ),u^{\delta
,\varepsilon ,\alpha }(\cdot ))-J(y_{0}^{\delta ,\varepsilon },z^{\delta
,\varepsilon }(\cdot ),u^{\delta ,\varepsilon }(\cdot ))\right] +\alpha
\varepsilon ^{\frac{1}{3}}\theta _{0}^{\delta ,\varepsilon ,\alpha }\sqrt{%
\left\vert y_{0}\right\vert ^{2}+2},
\end{align*}%
since the triangle inequality
\begin{eqnarray*}
&&d_{\mathcal{A}}(y_{0}^{\delta ,\varepsilon ,\alpha },z^{\delta
,\varepsilon ,\alpha }(\cdot ),u^{\delta ,\varepsilon ,\alpha }(\cdot ),%
\tilde{\theta}^{\varepsilon }(\cdot ))-d_{\mathcal{A}}(y_{0}^{\delta
,\varepsilon },z^{\delta ,\varepsilon }(\cdot ),u^{\delta ,\varepsilon
}(\cdot ),\tilde{\theta}^{\varepsilon }(\cdot )) \\
&\leq &d_{\mathcal{A}}\left( \left( y_{0}^{\delta ,\varepsilon ,\alpha
},z^{\delta ,\varepsilon ,\alpha }(\cdot ),u^{\delta ,\varepsilon ,\alpha
}(\cdot )\right) ,\left( y_{0}^{\delta ,\varepsilon },z^{\delta ,\varepsilon
}(\cdot ),u^{\delta ,\varepsilon }(\cdot )\right) \right) \leq \alpha \sqrt{%
\left\vert y_{0}\right\vert ^{2}+2}.
\end{eqnarray*}%
Note that
\begin{align*}
& \;J(y_{0}^{\delta ,\varepsilon ,\alpha },z^{\delta ,\varepsilon ,\alpha
}(\cdot ),u^{\delta ,\varepsilon ,\alpha }(\cdot ))-J(y_{0}^{\delta
,\varepsilon },z^{\delta ,\varepsilon }(\cdot ),u^{\delta ,\varepsilon
}(\cdot )) \\
=& \;\mathbb{E}\left[ \int_{0}^{1}l(t,X^{\delta ,\varepsilon ,\alpha
}(t),u^{\delta ,\varepsilon ,\alpha }(t))-l(t,X^{\delta ,\varepsilon
}(t),u^{\delta ,\varepsilon }(t))\mbox{\rm d}t\right] \\
& +\mathbb{E}\left[ \Xi (X^{\delta ,\varepsilon ,\alpha }(0),X^{\delta
,\varepsilon ,\alpha }(1))-\Xi (X^{\delta ,\varepsilon }(0),X^{\delta
,\varepsilon }(1))\right] \\
=& \;I_{1}+I_{2}.
\end{align*}%
We deal with $I_{1},$ $I_{2},$ respectively.%
\begin{eqnarray*}
I_{1} &=&\;\mathbb{E}\bigg[\int_{0}^{1}\bigg(l(t,X^{\delta ,\varepsilon
,\alpha }(t),u^{\delta ,\varepsilon ,\alpha }(t))-l(t,X^{\delta ,\varepsilon
}(t),u^{\delta ,\varepsilon }(t))\bigg)\mbox{\rm d}t \\
&&+\int_{0}^{1}l_{X}(t,X^{\delta ,\varepsilon }(t),u^{\delta ,\varepsilon
}(t))(X_{1}^{\delta ,\varepsilon ,\alpha }(t)+X_{2}^{\delta ,\varepsilon
,\alpha }(t))\mbox{\rm d}t \\
&&+\int_{0}^{1}\frac{1}{2}l_{XX}(t,X^{\delta ,\varepsilon }(t),u^{\delta
,\varepsilon }(t))(X_{1}^{\delta ,\varepsilon ,\alpha }(t))^{2}\mbox{\rm d}t
\\
&&+\int_{0}^{1}(l_{X}(t,X^{\delta ,\varepsilon ,\alpha }(t),u^{\delta
,\varepsilon ,\alpha }(t))-l_{X}(t,X^{\delta ,\varepsilon }(t),u^{\delta
,\varepsilon }(t)))(X^{\delta ,\varepsilon ,\alpha }(t)-X^{\delta
,\varepsilon }(t))\mbox{\rm d}t \\
&&+\int_{0}^{1}(l_{X}(t,X^{\delta ,\varepsilon }(t),u^{\delta ,\varepsilon
}(t)))(X^{\delta ,\varepsilon ,\alpha }(t)-X^{\delta ,\varepsilon
}(t)-X_{1}^{\delta ,\varepsilon ,\alpha }(t)-X_{2}^{\delta ,\varepsilon
,\alpha }(t))\mbox{\rm d}t\bigg] \\
&&+\mathbb{E}\bigg[\int_{0}^{1}\bigg (\beta \left[ l_{XX}(t,\beta X^{\delta
,\varepsilon }(t)+(1-\beta )X^{\delta ,\varepsilon ,\alpha }(t),u^{\delta
,\varepsilon ,\alpha }(t))\right. \\
&&\left. -l_{XX}(t,X^{\delta ,\varepsilon }(t),u^{\delta ,\varepsilon
,\alpha }(t))\right] (X^{\delta ,\varepsilon ,\alpha }(t)-X^{\delta
,\varepsilon }(t))^{2}\bigg )\mbox{\rm d}\beta \bigg]
\end{eqnarray*}%
\begin{eqnarray*}
&&+\frac{1}{2}\mathbb{E}\left[ \int_{0}^{1}(l_{XX}(t,X^{\delta ,\varepsilon
,\alpha }(t),u^{\delta ,\varepsilon ,\alpha }(t))-l_{XX}(t,X^{\delta
,\varepsilon }(t),u^{\delta ,\varepsilon }(t)))(X^{\delta ,\varepsilon
,\alpha }(t)-X^{\delta ,\varepsilon }(t))^{2}\right] \mbox{\rm d}t \\
&&+\frac{1}{2}\mathbb{E}\bigg[\int_{0}^{1}\bigg (l_{XX}(t,X^{\delta
,\varepsilon }(t),u^{\delta ,\varepsilon }(t))(X^{\delta ,\varepsilon
,\alpha }(t)-X^{\delta ,\varepsilon }(t)-X_{1}^{\delta ,\varepsilon ,\alpha
}(t)) \\
&&\times (X^{\delta ,\varepsilon ,\alpha }(t)-X^{\delta ,\varepsilon
}(t)+X_{1}^{\delta ,\varepsilon ,\alpha }(t))\bigg )\mbox{\rm d}t\bigg],
\end{eqnarray*}%
and
\begin{align*}
I_{2}=& \;\mathbb{E}\left[ \Xi (X^{\delta ,\varepsilon ,\alpha
}(0),X^{\delta ,\varepsilon ,\alpha }(1))-\Xi (X^{\delta ,\varepsilon
}(0),X^{\delta ,\varepsilon }(1))\right] \\
=& \;\mathbb{E}\left[ \Xi _{X(0)}(X^{\delta ,\varepsilon }(0),X^{\delta
,\varepsilon }(1))(X^{\delta ,\varepsilon ,\alpha }(0)-X^{\delta
,\varepsilon }(0))\right] \\
& +\mathbb{E}\left[ \Xi _{X(1)}(X^{\delta ,\varepsilon }(0),X^{\delta
,\varepsilon }(1))(X^{\delta ,\varepsilon ,\alpha }(1)-X^{\delta
,\varepsilon }(1))\right] \\
& +\mathbb{E}\left[ \frac{1}{2}\Xi _{X(0)X(0)}(X^{\delta ,\varepsilon
,\alpha }(0)-X^{\delta ,\varepsilon }(0))^{2}\right] +\mathbb{E}\left[ \frac{%
1}{2}\Xi _{X(1)X(1)}(X^{\delta ,\varepsilon ,\alpha }(1)-X^{\delta
,\varepsilon }(1))^{2}\right] \\
& +\mathbb{E}\left\langle D^{2}\Xi ^{\delta ,\varepsilon ,\alpha }\left(
\begin{array}{c}
X^{\delta ,\varepsilon ,\alpha }(0)-X^{\delta ,\varepsilon }(0) \\
X^{\delta ,\varepsilon ,\alpha }(1)-X^{\delta ,\varepsilon }(1)%
\end{array}%
\right) ,\left(
\begin{array}{c}
X^{\delta ,\varepsilon ,\alpha }(0)-X^{\delta ,\varepsilon }(0) \\
X^{\delta ,\varepsilon ,\alpha }(1)-X^{\delta ,\varepsilon }(1)%
\end{array}%
\right) \right\rangle ,
\end{align*}%
where
\begin{align*}
D^{2}\Xi ^{\delta ,\varepsilon ,\alpha }=& \int_{0}^{1}[\beta D^{2}\Xi
(\beta X^{\delta ,\varepsilon }(0)+(1-\beta )X^{\delta ,\varepsilon ,\alpha
}(0),\beta X^{\delta ,\varepsilon }(1)+(1-\beta X^{\delta ,\varepsilon
,\alpha }(1))) \\
& -D^{2}\Xi (X^{\delta ,\varepsilon }(0),X^{\delta ,\varepsilon }(1))]%
\mbox{\rm d}\beta .
\end{align*}%
Besides,
\begin{align*}
& \;\mathbb{E}\left\langle \binom{0}{\theta _{1}^{\delta ,\varepsilon
,\alpha }},\Pi (X^{\delta ,\varepsilon ,\alpha }(0),X^{\delta ,\varepsilon
,\alpha }(1))-\Pi (X^{\delta ,\varepsilon }(0),X^{\delta ,\varepsilon
}(1))\right\rangle \\
=& \;\mathbb{E}\left\langle \binom{0}{\theta _{1}^{\delta ,\varepsilon
,\alpha }},\Pi (0,X^{\delta ,\varepsilon ,\alpha }(1))-\Pi (0,X^{\delta
,\varepsilon }(1))\right\rangle \\
=& \;\mathbb{E}\left\langle \Pi _{X(1)}(0,X^{\delta ,\varepsilon }(1))\binom{%
0}{\theta _{1}^{\delta ,\varepsilon ,\alpha }}(X^{\delta ,\varepsilon
,\alpha }(1)-X^{\delta ,\varepsilon }(1))\right\rangle \\
& +\frac{1}{2}\mathbb{E}\left[ \left\langle \Pi _{X(1)X(1)}(0,X^{\delta
,\varepsilon }(1))\binom{0}{\theta _{1}^{\delta ,\varepsilon ,\alpha }}%
(X^{\delta ,\varepsilon ,\alpha }(1)-X^{\delta ,\varepsilon
}(1))^{2}\right\rangle \right] \\
& +\mathbb{E}\left[ D^{2}\Pi ^{\delta ,\varepsilon ,\alpha }\binom{0}{\theta
_{1}^{\delta ,\varepsilon ,\alpha }}(X^{\delta ,\varepsilon ,\alpha
}(1)-X^{\delta ,\varepsilon }(1))^{2}\right] ,
\end{align*}%
where
\begin{equation*}
D^{2}\Pi ^{\delta ,\varepsilon ,\alpha }=\int_{0}^{1}\beta \left[ \Pi
_{X(1)X(1)}(0,\beta X^{\delta ,\varepsilon }(1)+(1-\beta )X^{\delta
,\varepsilon ,\alpha }(1))-\Pi _{X(1)X(1)}(0,X^{\delta ,\varepsilon }(1))%
\right] \mbox{\rm d}\beta .
\end{equation*}%
Clearly, under assumptions (H1), we have%
\begin{eqnarray*}
&&-\alpha \sqrt{\left\vert y_{0}\right\vert ^{2}+2}\left( \sqrt{\delta }%
+\varepsilon ^{\frac{1}{3}}\theta _{0}^{\delta ,\varepsilon ,\alpha }\right)
\\
&\leq &\;\theta _{0}^{\delta ,\varepsilon ,\alpha }\left[ J(y_{0}^{\delta
,\varepsilon ,\alpha },z^{\delta ,\varepsilon ,\alpha }(\cdot ),u^{\delta
,\varepsilon ,\alpha }(\cdot ))-J(y_{0}^{\delta ,\varepsilon },z^{\delta
,\varepsilon }(\cdot ),u^{\delta ,\varepsilon }(\cdot ))\right] \\
&&+\mathbb{E}\left\langle \binom{0}{\theta _{1}^{\delta ,\varepsilon ,\alpha
}},\Pi (X^{\delta ,\varepsilon ,\alpha }(0),X^{\delta ,\varepsilon ,\alpha
}(1))-\Pi (X^{\delta ,\varepsilon }(0),X^{\delta ,\varepsilon
}(1))\right\rangle \\
&=&\;\theta _{0}^{\delta ,\varepsilon ,\alpha }\mathbb{E}\bigg[\int_{0}^{1}%
\bigg(l_{X}(t,X^{\delta ,\varepsilon }(t),u^{\delta ,\varepsilon
}(t))(X_{1}^{\delta ,\varepsilon ,\alpha }(t)+X_{2}^{\delta ,\varepsilon
,\alpha }(t)) \\
&&+l(t,X^{\delta ,\varepsilon ,\alpha }(t),u^{\delta ,\varepsilon ,\alpha
}(t))-l(t,X^{\delta ,\varepsilon }(t),u^{\delta ,\varepsilon }(t))+\frac{1}{2%
}l_{XX}(t,X^{\delta ,\varepsilon }(t),u^{\delta ,\varepsilon
}(t))(X_{1}^{\delta ,\varepsilon ,\alpha }(t))^{2}\bigg)\mbox{\rm d}t\bigg]
\\
&&+\mathbb{E}\left[ \Xi _{X(0)}(X^{\delta ,\varepsilon ,\alpha
}(0)-X^{\delta ,\varepsilon }(0))+\Xi _{X(1)}(X^{\delta ,\varepsilon ,\alpha
}(1)-X^{\delta ,\varepsilon }(1))\right] \\
&&+\mathbb{E}\left[ \frac{1}{2}\Xi _{X(0)X(0)}(X^{\delta ,\varepsilon
}(t),u^{\delta ,\varepsilon }(\cdot ))(X^{\delta ,\varepsilon ,\alpha
}(0)-X^{\delta ,\varepsilon }(0))^{2}\right] \\
&&+\mathbb{E}\left[ \frac{1}{2}\Xi _{X(1)X(1)}(X^{\delta ,\varepsilon
}(t),u^{\delta ,\varepsilon }(\cdot ))(X^{\delta ,\varepsilon ,\alpha
}(1)-X^{\delta ,\varepsilon }(1))^{2}\right] \\
&&+\mathbb{E}\left\langle \Pi _{X(1)}(0,X^{\delta ,\varepsilon }(1))\binom{0%
}{\theta _{1}^{\delta ,\varepsilon ,\alpha }}(X^{\delta ,\varepsilon ,\alpha
}(1)-X^{\delta ,\varepsilon }(1))\right\rangle \\
&&+\mathbb{E}\left[ \frac{1}{2}\left\langle \Pi _{X(1)X(1)}(0,X^{\delta
,\varepsilon }(1))\binom{0}{\theta _{1}^{\delta ,\varepsilon ,\alpha }}%
(X^{\delta ,\varepsilon ,\alpha }(1)-X^{\delta ,\varepsilon
}(1))^{2}\right\rangle \right] +o(\alpha ) \\
&=&\;\mathbb{E}\Bigg \{\int_{0}^{1}\bigg \{\theta _{0}^{\delta ,\varepsilon
,\alpha }\left[ l(t,X^{\delta ,\varepsilon ,\alpha }(t),u^{\delta
,\varepsilon ,\alpha }(t))-l(t,X^{\delta ,\varepsilon }(t),u^{\delta
,\varepsilon }(t))\right] \\
&&+\theta _{0}^{\delta ,\varepsilon ,\alpha }\left[ l_{X}(t,X^{\delta
,\varepsilon }(t),u^{\delta ,\varepsilon }(t))(X_{1}^{\delta ,\varepsilon
,\alpha }(t)+X_{2}^{\delta ,\varepsilon ,\alpha }(t))\right] \\
&&+\frac{1}{2}\theta _{0}^{\delta ,\varepsilon ,\alpha }l_{XX}(t,X^{\delta
,\varepsilon }(t),u^{\delta ,\varepsilon }(t))\left( X_{1}^{\delta
,\varepsilon ,\alpha }\left( t\right) \right) ^{2}\bigg \}\mbox{\rm d}t \\
&&+\sqrt{\alpha }\left\langle \theta _{0}^{\delta ,\varepsilon ,\alpha }\Xi
_{X(0)}\left( X^{\delta ,\varepsilon }\left( 0\right) ,X^{\delta
,\varepsilon }\left( 1\right) \right) ,\binom{0}{y_{0}}\right\rangle \\
&&+\frac{\alpha }{2}\left\langle \theta _{0}^{\delta ,\varepsilon ,\alpha
}\Xi _{X(0)X(0}\left( X^{\delta ,\varepsilon }\left( 0\right) ,X^{\delta
,\varepsilon }\left( 1\right) \right) \binom{0}{y_{0}},\binom{0}{y_{0}}%
\right\rangle \\
&&+\left( \theta _{0}^{\delta ,\varepsilon ,\alpha }\Xi _{X(1)}(X^{\delta
,\varepsilon }(0),X^{\delta ,\varepsilon }(1))+\Pi _{X(1)}(0,X^{\delta
,\varepsilon }(1))\binom{0}{\theta _{1}^{\delta ,\varepsilon ,\alpha }}%
\right) (X_{1}^{\delta ,\varepsilon ,\alpha }(1)+X_{2}^{\delta ,\varepsilon
,\alpha }(1)) \\
&&+\frac{1}{2}\Bigg[\theta _{0}^{\delta ,\varepsilon ,\alpha }\Xi
_{X(1)X(1)}(X^{\delta ,\varepsilon }(1),u^{\delta ,\varepsilon }(\cdot
))+\Pi _{X(1)X(1)}(0,X^{\delta ,\varepsilon }(1))\binom{0}{\theta
_{1}^{\delta ,\varepsilon ,\alpha }}\Bigg](X_{1}^{\delta ,\varepsilon
,\alpha }(1))^{2}\Bigg \}+o(\alpha ).
\end{eqnarray*}%
\begin{equation}  \label{6.45}
\end{equation}

\noindent Let us introduce the following the first order BSDEs:
\begin{equation*}
\left\{
\begin{array}{rcl}
-\mbox{\rm d}\tilde{\Phi}^{\delta ,\varepsilon ,\alpha }(t) & = & \left[
\mathcal{B}_{X}^{\delta ,\varepsilon }(t,\cdot )\tilde{\Phi}^{\delta
,\varepsilon ,\alpha }(t)+\Sigma _{X}^{\delta ,\varepsilon }(t,\cdot )\tilde{%
\Psi}^{\delta ,\varepsilon ,\alpha }(t)+\theta _{0}^{\delta ,\varepsilon
,\alpha }l_{X}^{\delta ,\varepsilon }(t,\cdot )\right] \mbox{\rm d}t-\tilde{%
\Psi}^{\delta ,\varepsilon ,\alpha }(t)\mbox{\rm d}W(t), \\[2mm]
\tilde{\Phi}^{\delta ,\varepsilon ,\alpha }(1) & = & \displaystyle\left[
\theta _{0}^{\delta ,\varepsilon ,\alpha }\Xi _{X(1)}(X^{\delta ,\varepsilon
}(0),X^{\delta ,\varepsilon }(1))+\Pi _{X(1)}(0,X^{\delta ,\varepsilon }(1))%
\binom{0}{\theta _{1}^{\delta ,\varepsilon ,\alpha }}\right] ,%
\end{array}%
\right.
\end{equation*}%
where $l_{X}^{\delta ,\varepsilon }(t,\cdot )=l_{X}(t,X^{\delta ,\varepsilon
}(t),u^{\delta ,\varepsilon }(t))$.

The second order BSDEs:
\begin{equation*}
\left\{
\begin{array}{rcl}
-\mbox{\rm d}\tilde{P}^{\delta ,\varepsilon ,\alpha }(t) & = & \left[
\mathcal{B}_{X}^{\delta ,\varepsilon }(t,\cdot )^{\top }\tilde{P}^{\delta
,\varepsilon ,\alpha }(t)+\tilde{P}^{\delta ,\varepsilon ,\alpha }(t)%
\mathcal{B}_{X}^{\delta ,\varepsilon }(t,\cdot )+\Sigma _{X}^{\delta
,\varepsilon }(t,\cdot )^{\top }\tilde{P}^{\delta ,\varepsilon ,\alpha
}(t)\Sigma _{X}^{\delta ,\varepsilon }(t,\cdot )\right. \\[2mm]
&  & \left. \Sigma _{X}^{\delta ,\varepsilon }(t,\cdot )^{\top }\tilde{Q}%
^{\delta ,\varepsilon ,\alpha }(t)+\tilde{Q}^{\delta ,\varepsilon ,\alpha
}(t)\Sigma _{X}^{\delta ,\varepsilon }(t,\cdot )+H_{XX}^{\delta ,\varepsilon
,\alpha }\right] \mbox{\rm d}t-\tilde{Q}^{\delta ,\varepsilon ,\alpha }(t)%
\text{d}W(t), \\[2mm]
\tilde{P}^{\delta ,\varepsilon ,\alpha }(1) & = & \displaystyle\left[ \theta
_{0}^{\delta ,\varepsilon ,\alpha }\Xi _{X(1)X(1)}(X^{\delta ,\varepsilon
}(0),X^{\delta ,\varepsilon }(1))+\Pi _{X(1)X(1)}(0,X^{\delta ,\varepsilon
}(1))\binom{0}{\theta _{1}^{\delta ,\varepsilon ,\alpha }}\right] ,%
\end{array}%
\right.
\end{equation*}%
where
\begin{equation*}
\tilde{H}_{XX}^{\delta ,\varepsilon ,\alpha }(t)=\tilde{H}_{XX}(t,\theta
_{0}^{\delta ,\varepsilon ,\alpha },X^{\delta ,\varepsilon ,\alpha
}(t),u^{\delta ,\varepsilon ,\alpha }(t),\tilde{\Phi}^{\delta ,\varepsilon
,\alpha }(t),\tilde{\Psi}^{\delta ,\varepsilon ,\alpha }(t)),
\end{equation*}%
with $\tilde{H}_{XX}(t,\theta _{0}^{\delta ,\varepsilon ,\alpha },X^{\delta
,\varepsilon }(t),u^{\delta ,\varepsilon }(t),\tilde{\Phi}^{\delta
,\varepsilon ,\alpha }(t),\tilde{\Psi}^{\delta ,\varepsilon ,\alpha }(t))$
is defined as follows:
\begin{equation*}
\tilde{H}(t,\theta ,X,v,p,k)=\left\langle p,\mathcal{B}(t,X,v)\right\rangle
+\left\langle k,\Sigma (t,X,v)\right\rangle +\theta l(t,X,v).
\end{equation*}%
Set $\mathcal{Y}^{\delta ,\varepsilon ,\alpha }(\cdot )=X_{1}^{\delta
,\varepsilon ,\alpha }(\cdot )X_{1}^{\delta ,\varepsilon ,\alpha }(\cdot ).$
Then,
\begin{equation*}
\left\{
\begin{array}{rcl}
\mbox{\rm d}\mathcal{Y}^{\delta ,\varepsilon ,\alpha }(t) & = & \Bigg \{%
\mathcal{B}_{X}^{\delta }(t,\cdot )\mathcal{Y}^{\delta ,\varepsilon ,\alpha
}(t)+\mathcal{Y}^{\delta ,\varepsilon ,\alpha }(t)\mathcal{B}_{X}^{\delta
}(t,\cdot )^{\top }+\Sigma _{X}^{\delta }(t,\cdot )\mathcal{Y}^{\delta
,\varepsilon ,\alpha }(t)\Sigma _{X}^{\delta }(t,\cdot )^{\top } \\[2mm]
&  & +\bigg[\triangle \Sigma ^{\delta }(t,\cdot )\triangle \Sigma ^{\delta
}(t,\cdot )^{\top }+\Sigma _{X}^{\delta }(t,\cdot )X_{1}^{\delta
,\varepsilon ,\alpha }(t)\triangle \Sigma ^{\delta }(t,\cdot )^{\top } \\%
[2mm]
&  & +\triangle \Sigma ^{\delta }(t,\cdot )X_{1}^{\delta ,\varepsilon
,\alpha }(t)\Sigma ^{\delta }(t,\cdot )^{\top }\bigg]I_{S_{\alpha }}(t)\Bigg
\}\mbox{\rm d}t \\[2mm]
&  & +\bigg \{\Sigma _{X}^{\delta }(t,\cdot )\mathcal{Y}^{\delta
,\varepsilon ,\alpha }(t)+\mathcal{Y}^{\delta ,\varepsilon ,\alpha
}(t)\Sigma _{X}^{\delta }(t,\cdot )^{\top } \\[2mm]
&  & +\left[ X_{1}^{\delta ,\varepsilon ,\alpha }(t)\triangle \Sigma
^{\delta }(t,\cdot )^{\top }+\triangle \Sigma ^{\delta }(t,\cdot
)X_{1}^{\delta ,\varepsilon ,\alpha }(t)^{\top }\right] I_{S_{\alpha }}(t)%
\bigg \}\mbox{\rm d}W(t), \\[2mm]
\mathcal{Y}^{\delta ,\varepsilon ,\alpha }(0) & = & \left(
\begin{array}{cc}
0 & 0 \\
0 & \alpha y_{0}^{2}%
\end{array}%
\right) ,%
\end{array}%
\right.
\end{equation*}%
Applying Itô's formula to $\left\langle \tilde{\Phi}^{\delta ,\varepsilon
,\alpha }(\cdot ),X_{1}^{\delta ,\varepsilon ,\alpha }(\cdot )+X_{2}^{\delta
,\varepsilon ,\alpha }(\cdot )\right\rangle $ and $P^{\delta ,\varepsilon
,\alpha }(\cdot )\mathcal{Y}^{\delta ,\varepsilon ,\alpha }(\cdot )$
respectively, we have
\begin{align}
& \;\mathbb{E}\left[ \left\langle \tilde{\Phi}^{\delta ,\varepsilon ,\alpha
}(1),X_{1}^{\delta ,\varepsilon ,\alpha }(1)+X_{2}^{\delta ,\varepsilon
,\alpha }(1)\right\rangle \right] -\mathbb{E}\left[ \left\langle \tilde{\Phi}%
^{\delta ,\varepsilon ,\alpha }(0),\binom{0}{\sqrt{\alpha }y_{0}}%
\right\rangle \right]  \notag  \label{6.49} \\
=& \;\mathbb{E}\bigg[\int_{0}^{1}-\left\langle \theta _{0}^{\delta
,\varepsilon ,\alpha }l_{X}^{\delta ,\alpha }(t,\cdot ),(X_{2}^{\delta
,\varepsilon ,\alpha }(t)+X_{1}^{\delta ,\varepsilon ,\alpha
}(t))\right\rangle  \notag \\
& +\left\langle \tilde{\Phi}^{\delta ,\varepsilon ,\alpha }(t),\triangle
\mathcal{B}^{\delta ,\varepsilon }(t,\cdot )I_{S_{\alpha }}(t)+\frac{1}{2}%
\mathcal{B}_{XX}^{\delta ,\varepsilon }(t,\cdot )(X_{1}^{\delta ,\varepsilon
,\alpha }(t))^{2}\right\rangle  \notag \\
& +\left\langle \tilde{\Psi}^{\delta ,\varepsilon ,\alpha }(t),\triangle
\Sigma ^{\delta ,\varepsilon }(t,\cdot )X_{1}^{\delta ,\varepsilon ,\alpha
}(t)I_{S_{\alpha }}(t)+\frac{1}{2}\Sigma _{XX}^{\delta ,\varepsilon
}(t,\cdot )(X_{1}^{\delta ,\varepsilon ,\alpha }(t))^{2}\right\rangle \text{d%
}t\bigg]+o(\alpha ).
\end{align}%
and
\begin{align}
& \;\mathbb{E}\left[ \text{tr}\left[ P^{\delta ,\varepsilon ,\alpha }(1)%
\mathcal{Y}^{\delta ,\varepsilon ,\alpha }(1)\right] -\left\langle P^{\delta
,\varepsilon ,\alpha }(1)\binom{0}{\sqrt{\varepsilon }y_{0}},\binom{0}{\sqrt{%
\varepsilon }y_{0}}\right\rangle \right]  \notag  \label{6.50} \\
=& \;\mathbb{E}\left\{ \int_{0}^{1}\text{tr}\left[ \triangle \Sigma ^{\delta
,\varepsilon }(t,\cdot )^{\top }P^{\delta ,\varepsilon ,\alpha }(t)\triangle
\Sigma ^{\delta ,\varepsilon }(t,\cdot )-\left\langle H_{XX}^{\delta
,\varepsilon ,\alpha }(t)X_{1}^{\delta ,\varepsilon ,\alpha
}(t),X_{1}^{\delta ,\varepsilon ,\alpha }(t)\right\rangle \right]
\mbox{\rm
d}t\right\} +o(\alpha ).
\end{align}%
Then, from (\ref{6.45}), (\ref{6.49}) and (\ref{6.50}), we obtain%
\begin{eqnarray}
&&-\alpha \sqrt{\left\vert y_{0}\right\vert ^{2}+2}\left( \sqrt{\delta }%
+\varepsilon ^{\frac{1}{3}}\theta _{0}^{\delta ,\varepsilon ,\alpha }\right)
\notag \\
&\leq &\;\mathbb{E}\bigg[\int_{0}^{1}\bigg [\theta _{0}^{\delta ,\varepsilon
,\alpha }\left[ l(t,X^{\delta ,\varepsilon ,\alpha }(t),u^{\delta
,\varepsilon ,\alpha }(\cdot ))-l(t,X^{\delta ,\varepsilon }(t),u^{\delta
,\varepsilon }(\cdot ))\right]  \notag \\
&&+\left\langle \tilde{\Phi}^{\delta ,\varepsilon ,\alpha }(t),\triangle
\mathcal{B}^{\delta ,\varepsilon }(t,\cdot )\right\rangle +\left\langle
\tilde{\Psi}^{\delta ,\varepsilon ,\alpha }(t),\triangle \Sigma ^{\delta
,\varepsilon }(t,\cdot )\right\rangle  \notag \\
&&+\frac{1}{2}\triangle \Sigma ^{\delta ,\varepsilon }(t,\cdot )^{\top
}P^{\delta ,\varepsilon ,\alpha }(t)\triangle \Sigma ^{\delta ,\varepsilon
}(t,\cdot )\bigg ]\mbox{\rm d}t\bigg]  \notag \\
&&+\mathbb{E}\left[ \sqrt{\alpha }\left\langle \theta _{0}^{\delta
,\varepsilon ,\alpha }\Xi _{X(0)}\left( X^{\delta ,\varepsilon }\left(
0\right) ,X^{\delta ,\varepsilon }\left( 1\right) \right) +\tilde{\Phi}%
^{\delta ,\varepsilon ,\alpha }(0),\binom{0}{y_{0}}\right\rangle \right]
\notag \\
&&+\mathbb{E}\left[ \frac{\alpha }{2}\left\langle \left( \theta _{0}^{\delta
,\varepsilon ,\alpha }\Xi _{X(0)X(0}\left( X^{\delta ,\varepsilon }\left(
0\right) ,X^{\delta ,\varepsilon }\left( 1\right) \right) +P^{\delta
,\varepsilon ,\alpha }(0)\right) \binom{0}{y_{0}},\binom{0}{y_{0}}%
\right\rangle \right] +o(\alpha ).  \label{6.51}
\end{eqnarray}%
To derive the adjoint equations, in (\ref{6.51})$,$ dividing $\sqrt{\alpha }$
and then sending $\alpha \rightarrow 0,$ followed by sending $\delta
\rightarrow 0,$ we get%
\begin{align}
0\leq & \;\mathbb{E}\left\langle \theta _{0}^{\varepsilon }\Xi _{X(0)}(%
\tilde{X}^{\varepsilon }(0),\tilde{X}^{\varepsilon }(1))+\tilde{\Phi}%
^{\varepsilon }(0),\binom{0}{y_{0}}\right\rangle  \notag \\
=& \;\mathbb{E}\left\langle \theta _{0}^{\varepsilon }\binom{0}{\gamma _{y}(%
\tilde{y}_{0}^{\varepsilon })}+\tilde{\Phi}^{\varepsilon }(0),\binom{0}{y_{0}%
}\right\rangle .  \label{6.52}
\end{align}%
From continuous dependence of the solution of BSDEs on parameters $\left(
\theta _{0}^{\delta ,\varepsilon ,\alpha },\theta _{1}^{\delta ,\varepsilon
,\alpha }\right) $, we get
\begin{eqnarray*}
\left( \theta _{0}^{\delta ,\varepsilon ,\alpha },\theta _{1}^{\delta
,\varepsilon ,\alpha }\right) &\rightarrow &\left( \theta _{0}^{\varepsilon
},\theta _{1}^{\varepsilon }\right) \in \mathbb{R}\times L_{\mathcal{F}%
_{1}}^{2}(\Omega ;\mathbb{R}),\text{ weakly,} \\
\left( \tilde{\Phi}^{\delta ,\varepsilon ,\alpha }(\cdot ),\tilde{\Psi}%
^{\delta ,\varepsilon ,\alpha }(\cdot )\right) &\rightarrow &\left( \tilde{%
\Phi}^{\varepsilon }(\cdot ),\tilde{\Psi}^{\varepsilon }(\cdot )\right) ,%
\text{ in }\mathcal{M}^{2}(0,1;\mathbb{R}), \\
\left( \tilde{P}^{\delta ,\varepsilon ,\alpha }(\cdot ),\tilde{Q}^{\delta
,\varepsilon ,\alpha }(\cdot )\right) &\rightarrow &\left( \tilde{P}%
^{\varepsilon }(\cdot ),\tilde{Q}^{\varepsilon }(\cdot )\right) ,\text{ in }%
\mathcal{M}^{2}(0,1;\mathbb{R}),\text{ as }\delta \rightarrow 0,\alpha
\rightarrow 0.
\end{eqnarray*}%
Denote
\begin{equation*}
\tilde{\Phi}^{\varepsilon }(\cdot )=\binom{\tilde{p}^{\varepsilon }(\cdot )}{%
\tilde{q}^{\varepsilon }(\cdot )},\quad \tilde{\Psi}^{\varepsilon }(\cdot )=%
\binom{\tilde{k}^{\varepsilon }(\cdot )}{\tilde{h}^{\varepsilon }(\cdot )}.
\end{equation*}%
Then, from (\ref{6.52}), we derive that
\begin{equation}
\binom{\tilde{p}^{\varepsilon }(0)}{\tilde{q}^{\varepsilon }(0)}=\binom{0}{%
-\theta _{0}^{\varepsilon }\mathbb{E}\gamma _{y}(\tilde{y}^{\varepsilon }(1))%
}.\quad  \label{A18}
\end{equation}%
Note that%
\begin{equation*}
\left\{
\begin{array}{ccc}
\displaystyle\Xi _{X(1)}(\tilde{X}^{\varepsilon }(0),\tilde{X}^{\varepsilon
}(1)) & = & \displaystyle\binom{\phi _{x}(\tilde{x}^{\varepsilon }(1))}{0},
\\
\displaystyle\Xi _{X(1)X(1)}(\tilde{X}^{\varepsilon }(0),\tilde{X}%
^{\varepsilon }(1)) & = & \left(
\begin{array}{cc}
\phi _{xx}(\tilde{x}^{\varepsilon }(1)) & 0 \\
0 & 0%
\end{array}%
\right) ,%
\end{array}%
\right.
\end{equation*}%
and
\begin{equation*}
\left\{
\begin{array}{ccc}
\displaystyle\Pi _{X(1)}(0,\tilde{X}^{\varepsilon }(1))\binom{0}{\theta
_{1}^{\varepsilon }} & = & \displaystyle\binom{-M\theta _{1}^{\varepsilon }}{%
\theta _{1}^{\varepsilon }}, \\
\displaystyle\Pi _{X(1)X(1)}(0,\tilde{X}^{\varepsilon }(1))\binom{0}{\theta
_{1}^{\varepsilon }} & = & \left(
\begin{array}{cc}
0 & 0 \\
0 & 0%
\end{array}%
\right) .%
\end{array}%
\right.
\end{equation*}%
Next, for the first and second order BSDEs for Problem $(\tilde{C}%
^{\varepsilon })$, we have%
\begin{equation}
\binom{\tilde{p}^{\varepsilon }(1)}{\tilde{q}^{\varepsilon }(1)}=\binom{%
\theta _{0}^{\varepsilon }\phi _{x}(\tilde{x}^{\varepsilon }(1))-M\theta
_{1}^{\varepsilon }}{\theta _{1}^{\varepsilon }},  \label{A19}
\end{equation}%
and
\begin{align}
\tilde{P}^{\varepsilon }(1)=& \left[ \theta _{0}^{\varepsilon }\Xi
_{X(1)X(1)}(\tilde{X}^{\varepsilon }(0),\tilde{X}^{\varepsilon }(1))+\Pi
_{X(1)X(1)}(0,\tilde{X}^{\varepsilon }(1))\binom{0}{\theta _{1}^{\varepsilon
}}\right]  \notag \\
=& \left(
\begin{array}{cc}
\theta _{0}^{\varepsilon }\phi _{xx}(\tilde{x}^{\varepsilon }(1)) & 0 \\
0 & 0%
\end{array}%
\right) .  \label{A20}
\end{align}%
where the first and second order BSDEs are
\begin{equation*}
\left\{
\begin{array}{rcl}
-\mbox{\rm d}\tilde{\Phi}^{\varepsilon }(t) & = & \left[ \mathcal{B}%
_{X}(t,\cdot )\tilde{\Phi}^{\varepsilon }(t)+\Sigma _{X}(t,\cdot )\tilde{\Psi%
}^{\varepsilon }(t)+\theta _{0}^{\varepsilon }l_{X}(t,\cdot )\right]
\mbox{\rm
d}t \\[2mm]
&  & -\tilde{\Psi}^{\varepsilon }(t)\mbox{\rm d}W(t), \\[2mm]
\tilde{\Phi}^{\varepsilon }(1) & = & \displaystyle\left[ \theta
_{0}^{\varepsilon }\Xi _{X(1)}(\tilde{X}^{\varepsilon }(0),\tilde{X}%
^{\varepsilon }(1))+\Pi _{X(1)}(0,\tilde{X}^{\varepsilon }(1))\binom{0}{%
\theta _{1}^{\varepsilon }}\right] ,%
\end{array}%
\right.
\end{equation*}%
and
\begin{equation*}
\left\{
\begin{array}{rcl}
-\mbox{\rm d}\tilde{P}^{\varepsilon }(t) & = & \left[ \mathcal{B}%
_{X}(t,\cdot )^{\top }\tilde{P}^{\varepsilon }(t)+\tilde{P}^{\varepsilon }(t)%
\mathcal{B}_{X}(t,\cdot )+\Sigma _{X}(t,\cdot )^{\top }\tilde{P}%
^{\varepsilon }(t)\Sigma _{X}(t,\cdot )\right. \\[2mm]
&  & \left. \Sigma _{X}(t,\cdot )^{\top }\tilde{Q}^{\varepsilon }(t)+\tilde{Q%
}^{\varepsilon }(t)\Sigma _{X}(t,\cdot )+H_{XX}\right] \mbox{\rm d}t-\tilde{Q%
}^{\varepsilon }(t)\mbox{\rm d}W(t), \\[2mm]
\tilde{P}^{\varepsilon }(1) & = & \displaystyle\left[ \theta
_{0}^{\varepsilon }\Xi _{X(1)X(1)}(\tilde{X}^{\varepsilon }(0),\tilde{X}%
^{\varepsilon }(1))+\Pi _{X(1)X(1)}(0,\tilde{X}^{\varepsilon }(1))\binom{0}{%
\theta _{1}^{\varepsilon }}\right] .%
\end{array}%
\right.
\end{equation*}%
Then using a standard argument of \cite{YZ}, taking $y_{0}=0,$ we have the
following variational inequality:
\begin{align}
-\sqrt{2}\varepsilon ^{\frac{1}{3}}\theta _{0}^{\varepsilon }& \leq \theta
_{0}^{\varepsilon }\left[ l(t,\tilde{x}^{\varepsilon }(t),\tilde{y}%
^{\varepsilon }(t),u)-l(t,\tilde{x}^{\varepsilon }(t),\tilde{y}^{\varepsilon
}(t),\tilde{u}^{\varepsilon }(t))\right]  \notag \\
& +\left\langle \tilde{\Phi}^{\varepsilon }(t),\mathcal{B}(t,\tilde{x}%
^{\varepsilon }(t),\tilde{y}^{\varepsilon }(t),u)-\mathcal{B}(t,\tilde{x}%
^{\varepsilon }(t),\tilde{y}^{\varepsilon }(t),\tilde{u}^{\varepsilon
}(t))\right\rangle  \notag \\
& +\left\langle \tilde{\Psi}^{\varepsilon }(t),\Sigma (t,\tilde{x}%
^{\varepsilon }(t),\tilde{y}^{\varepsilon }(t),u,z)-\Sigma (t,\tilde{x}%
^{\varepsilon }(t),\tilde{y}^{\varepsilon }(t),\tilde{u}^{\varepsilon }(t),%
\tilde{z}^{\varepsilon }\left( t\right) )\right\rangle  \notag \\
& +\frac{1}{2}(\Sigma (t,\tilde{x}^{\varepsilon }(t),\tilde{y}^{\varepsilon
}(t),u,z)-\Sigma (t,\tilde{x}^{\varepsilon }(t),\tilde{y}^{\varepsilon }(t),%
\tilde{u}^{\varepsilon }(t),\tilde{z}^{\varepsilon }\left( t\right) ))^{\top
}  \notag \\
& \times \tilde{P}^{\varepsilon }(t)(\Sigma (t,\tilde{x}^{\varepsilon }(t),%
\tilde{y}^{\varepsilon }(t),u)-\Sigma (t,\tilde{x}^{\varepsilon }(t),\tilde{y%
}^{\varepsilon }(t),\tilde{u}^{\varepsilon }(t)))+o(\alpha ),\text{ }  \notag
\\
\forall u& \in \mathbb{U},\text{ }\forall z\in \mathbb{N},\text{ }u\in
\mathbb{U}\text{, .a.e, a.s..}  \label{6.53}
\end{align}%
Then (\ref{6.53}) can be rewrote as
\begin{eqnarray}
-\sqrt{2}\varepsilon ^{\frac{1}{3}}\theta _{0}^{\varepsilon } &<&\theta
_{0}^{\varepsilon }\left[ l(t,\tilde{x}^{\varepsilon }(t),\tilde{y}%
^{\varepsilon }(t),u)-l(t,\tilde{x}^{\varepsilon }(t),\tilde{y}^{\varepsilon
}(t),\tilde{u}^{\varepsilon }(t))\right]  \notag \\
&&+\left\langle \tilde{p}^{\varepsilon }(t),B(t)(u-\tilde{u}^{\varepsilon
}(t))\right\rangle +\left\langle \tilde{k}^{\varepsilon }(t),D(t)(u-\tilde{u}%
^{\varepsilon }(t))\right\rangle  \notag \\
&&+\left\langle \tilde{h}^{\varepsilon }(t),z-\tilde{z}^{\varepsilon }\left(
t\right) \right\rangle -\left\langle \tilde{q}^{\varepsilon }(t),c(t)(u-%
\tilde{u}^{\varepsilon }(t))\right\rangle +\frac{1}{2}D^{2}(t)(u-\tilde{u}%
^{\varepsilon }(t))^{2}\tilde{P}_{1}^{\varepsilon }(t)  \notag \\
&&+\frac{1}{2}\binom{D(t)(u-\tilde{u}^{\varepsilon }(t)}{z-\tilde{z}%
^{\varepsilon }\left( t\right) }^{\top }\tilde{P}^{\varepsilon }\left(
t\right) \binom{D(t)(u-\tilde{u}^{\varepsilon }(t)}{z-\tilde{z}^{\varepsilon
}\left( t\right) },\text{ }  \notag \\
\forall z &\in &\mathbb{N},\text{ }u\in \mathbb{U}\text{, a.e, a.s..}
\label{A22}
\end{eqnarray}%
Taking $u\left( t\right) =\tilde{u}^{\varepsilon }\left( t\right) ,$ $%
z\left( t\right) =\tilde{z}^{\varepsilon }\left( t\right) +\epsilon z_{0},$ $%
\forall z_{0}\in \mathbb{N}$, then dividing by sending $\epsilon \rightarrow
0,$ we have%
\begin{equation*}
-\sqrt{2}\varepsilon ^{\frac{1}{3}}\theta _{0}^{\varepsilon }\leq
\left\langle \tilde{h}^{\varepsilon }(t),z_{0}\right\rangle .
\end{equation*}%
Hence, we derive that $\tilde{h}^{\varepsilon }\left( t\right) \equiv 0$
since $\theta _{0}^{\varepsilon }\geq 0.$ From (\ref{A18})-(\ref{A20}) we
get
\begin{equation}
\left\{
\begin{array}{rcl}
-\mbox{\rm d}\tilde{p}^{\varepsilon }(t) & = & \left[ A(t)\tilde{p}%
^{\varepsilon }(t)-a(t)\tilde{q}^{\varepsilon }(t)+C(t)\tilde{k}%
^{\varepsilon }(t)+\theta _{0}^{\varepsilon }l_{x}(t,\cdot )\right] \text{d}%
t-\tilde{k}^{\varepsilon }(t)\mbox{\rm d}W(t), \\[2mm]
\text{ \ }\mbox{\rm d}\tilde{q}^{\varepsilon }(t) & = & \left[ -b(t)\tilde{q}%
^{\varepsilon }(t)-\theta _{0}^{\varepsilon }l_{y}(t,\cdot )\right] %
\mbox{\rm d}t, \\[2mm]
\text{ \ }\tilde{p}^{\varepsilon }(1) & = & \theta _{0}^{\varepsilon }\phi
_{x}(x^{\varepsilon }(1))-M\theta _{1}^{\varepsilon },\text{ }\tilde{q}%
(0)=-\theta _{0}^{\varepsilon }\mathbb{\gamma }_{y}(\tilde{y}^{\varepsilon
}(0)),\text{ }\tilde{q}^{\varepsilon }(1)=\theta _{1}^{\varepsilon },%
\end{array}%
\right.  \label{6.531}
\end{equation}%
and
\begin{equation}
\left\{
\begin{array}{rcl}
-\mbox{\rm d}\tilde{P}^{\varepsilon }(t) & = & \left[ \mathcal{B}%
_{X}(t,\cdot )^{\top }\tilde{P}^{\varepsilon }(t)+\tilde{P}^{\varepsilon }(t)%
\mathcal{B}_{X}(t,\cdot )+\Sigma _{X}(t,\cdot )^{\top }\tilde{P}%
^{\varepsilon }(t)\Sigma _{X}(t,\cdot )\right. \\[2mm]
&  & \left. +\Sigma _{X}(t,\cdot )^{\top }\tilde{Q}^{\varepsilon }(t)+\tilde{%
Q}^{\varepsilon }(t)\Sigma _{X}(t,\cdot )+H_{XX}(t,\cdot )\right]
\mbox{\rm
d}t-\tilde{Q}^{\varepsilon }(t)\mbox{\rm d}W(t), \\[2mm]
\text{ \ }\tilde{P}^{\varepsilon }(1) & = & \left(
\begin{array}{cc}
\theta _{0}^{\varepsilon }\phi _{xx}(\tilde{x}^{\varepsilon }(1)) & 0 \\
0 & 0%
\end{array}%
\right) ,%
\end{array}%
\right.  \label{6.532}
\end{equation}%
where
\begin{equation*}
H_{XX}(t,\cdot )=H_{XX}(t,\tilde{x}^{\varepsilon }(t),\tilde{y}^{\varepsilon
}(t),\tilde{u}^{\varepsilon }(t),\tilde{p}^{\varepsilon }(t),\tilde{q}%
^{\varepsilon }(t),\tilde{k}^{\varepsilon }(t),\theta _{0}^{\varepsilon }),
\end{equation*}%
and the Hamiltonian function $H:\left[ 0,T\right] \times \mathbb{R}\times
\mathbb{R}\times \mathbb{U}\times \mathbb{R}\times \mathbb{R}\times \mathbb{%
R\times R\rightarrow R}$ is defined as follows:%
\begin{align*}
H(t,x,y,u,p,q,k,\theta )\triangleq & \left\langle p,A(t)x+B(t)u\right\rangle
-\left\langle q,a(t)x+b(t)y+c(t)u\right\rangle \\
& +\left\langle k,C(t)x+D(t)u\right\rangle +\theta l(t,x,y,u).
\end{align*}%
Taking $y_{0}=0$ and $z\left( t\right) =\tilde{z}^{\varepsilon }(t)$ in (\ref%
{A22})$,$ we have the following variational inequality:%
\begin{eqnarray}
&&\left\langle \tilde{p}^{\varepsilon }(t),B(t)(u-\tilde{u}^{\varepsilon
}(t))\right\rangle +\left\langle \tilde{k}^{\varepsilon }(t),D(t)(u-\tilde{u}%
^{\varepsilon }(t))\right\rangle  \notag \\
&&-\left\langle \tilde{q}^{\varepsilon }(t),c(t)(u-\tilde{u}^{\varepsilon
}(t))\right\rangle  \notag \\
&&+\theta _{0}^{\varepsilon }\left[ l(t,\tilde{x}^{\varepsilon }(t),\tilde{y}%
^{\varepsilon }(t),u)-l(t,\tilde{x}^{\varepsilon }(t),\tilde{y}^{\varepsilon
}(t),\tilde{u}^{\varepsilon }(t))\right]  \notag \\
&&+\frac{1}{2}D^{2}(t)(u-\tilde{u}^{\varepsilon }(t))^{2}\tilde{P}%
_{1}^{\varepsilon }(t)  \notag \\
&\geq &-\sqrt{2}\varepsilon ^{\frac{1}{3}}\theta _{0}^{\varepsilon },\text{
a.e, a.s.}.  \label{6.54}
\end{eqnarray}%
Now consider (\ref{6.531})-(\ref{6.532}) again but only $(\tilde{x}%
^{\varepsilon }(\cdot ),\tilde{y}^{\varepsilon }(\cdot ),\tilde{z}%
^{\varepsilon }(\cdot ),\tilde{u}^{\varepsilon }(\cdot ))$ replaced by $%
(x^{\varepsilon }(\cdot ),y^{\varepsilon }(\cdot ),z^{\varepsilon }(\cdot
),u^{\varepsilon }(\cdot )).$ We need to derive an estimate for the term
similar to the right hand side of (\ref{6.54}) with all $(\tilde{x}%
^{\varepsilon }(\cdot ),\tilde{y}^{\varepsilon }(\cdot ),\tilde{z}%
^{\varepsilon }(\cdot ),\tilde{u}^{\varepsilon }(\cdot ))$ replaced by $%
(x^{\varepsilon }(\cdot ),y^{\varepsilon }(\cdot ),z^{\varepsilon }(\cdot
),u^{\varepsilon }(\cdot )).$ To this end, we first estimate the following
difference:%
\begin{align*}
& \;\mathbb{E}\left[ \int_{0}^{1}D(t)\left[ (u-\tilde{u}^{\varepsilon }(t))%
\tilde{k}^{\varepsilon }(t)\right] \mbox{\rm d}t\right] -\mathbb{E}\left[
\int_{0}^{1}D(t)\left[ (u-u^{\varepsilon }(t))k^{\varepsilon }(t)\right] %
\mbox{\rm d}t\right] \\
=& \;\mathbb{E}\left[ \int_{0}^{1}D(t)(u-\tilde{u}^{\varepsilon }(t))(\tilde{%
k}^{\varepsilon }(t)-k^{\varepsilon }(t))\mbox{\rm d}t\right] +\mathbb{E}%
\left[ \int_{0}^{1}D(t)(u^{\varepsilon }(t)-\tilde{u}^{\varepsilon
}(t))k^{\varepsilon }(t)\mbox{\rm d}t\right]
\end{align*}%
with
\begin{align*}
I_{1}=& \;\mathbb{E}\left[ \int_{0}^{1}D(t)(u-\tilde{u}^{\varepsilon }(t))(%
\tilde{k}^{\varepsilon }(t)-k^{\varepsilon }(t))\mbox{\rm d}t\right] , \\
I_{2}=& \;\mathbb{E}\left[ \int_{0}^{1}D(t)(u^{\varepsilon }(t)-\tilde{u}%
^{\varepsilon }(t))k^{\varepsilon }(t)\mbox{\rm d}t\right] .
\end{align*}%
Due to Lemma \ref{lem:16}, for any $1<\tau <2$ and $0<\beta <1$ satisfying $%
(1+\beta )\tau <2,$ there is a constant $C>0$ such that
\begin{align*}
I_{1}\leq & \;\left( \mathbb{E}\int_{0}^{1}\left\vert \tilde{k}^{\varepsilon
}(t)-k^{\varepsilon }(t)\right\vert ^{\tau }\mbox{\rm d}t\right) ^{\frac{1}{%
\tau }}\times \left( \mathbb{E}\int_{0}^{1}\left\vert u-\tilde{u}%
^{\varepsilon }(t)\right\vert ^{\frac{\tau }{\tau -1}}\mbox{\rm
d}t\right) ^{\frac{\tau -1}{\tau }} \\
\leq & \;C\left( d(u^{\varepsilon }(t)-\tilde{u}^{\varepsilon }(t))^{\frac{%
\tau \beta }{2}}\right) ^{\frac{1}{\tau }}\times \left( \mathbb{E}%
\int_{0}^{1}(\left\vert u\right\vert ^{\frac{\tau }{\tau -1}}+\left\vert
\tilde{u}^{\varepsilon }(t)\right\vert ^{\frac{\tau }{\tau -1}})\mbox{\rm
d}t\right) ^{\frac{\tau -1}{\tau }} \\
\leq & \;C\varepsilon ^{\frac{\beta }{3}},
\end{align*}%
and
\begin{eqnarray*}
I_{2} &\leq &\;C\left( \mathbb{E}\int_{0}^{1}\left\vert k^{\varepsilon
}(t)\right\vert ^{2}\mbox{\rm d}t\right) ^{\frac{1}{2}}\left( \mathbb{E}%
\int_{0}^{1}\left\vert u^{\varepsilon }(t)-\tilde{u}^{\varepsilon
}(t)\right\vert ^{^{2}}I_{u^{\varepsilon }(t)\neq \tilde{u}^{\varepsilon
}(t)}(t)\mbox{\rm d}t\right) ^{\frac{1}{2}} \\
&\leq &\;C\left( \mathbb{E}\int_{0}^{1}\left\vert u^{\varepsilon }(t)-\tilde{%
u}^{\varepsilon }(t)\right\vert ^{^{4}}\mbox{\rm d}t\right) ^{\frac{1}{4}%
}\left( \mathbb{E}\int_{0}^{1}I_{u^{\varepsilon }(t)\neq \tilde{u}%
^{\varepsilon }(t)}(t)\mbox{\rm d}t\right) ^{\frac{1}{4}} \\
&\leq &\;C\left( \mathbb{E}\int_{0}^{1}\left\vert u^{\varepsilon
}(t)\right\vert ^{4}+\left\vert \tilde{u}^{\varepsilon }(t)\right\vert
^{^{4}}\mbox{\rm d}t\right) ^{\frac{1}{4}}\left( \mathbb{E}%
\int_{0}^{1}I_{u^{\varepsilon }(t)\neq \tilde{u}^{\varepsilon }(t)}(t)%
\mbox{\rm d}t\right) ^{\frac{1}{4}} \\
&\leq &\;Cd(u^{\varepsilon }(\cdot ),\tilde{u}^{\varepsilon }(\cdot ))^{%
\frac{1}{4}} \\
&\leq &\;C\varepsilon ^{\frac{1}{6}} \\
&\leq &C\varepsilon ^{\frac{\beta }{3}}.
\end{eqnarray*}%
Similarly,
\begin{align*}
& \int_{0}^{1}\Big(\left\langle \tilde{p}^{\varepsilon }(t),B(t)(u-\tilde{u}%
^{\varepsilon }(t))\right\rangle -\left\langle p^{\varepsilon
}(t),B(t)(u-u^{\varepsilon }(t))\right\rangle \\
& +\left\langle q^{\varepsilon }(t),c(t)(u-u^{\varepsilon }(t))\right\rangle
-\left\langle \tilde{q}^{\varepsilon }(t),c(t)(u-\tilde{u}^{\varepsilon
}(t))\right\rangle \\
& +l(t,\tilde{x}^{\varepsilon }(t),\tilde{y}^{\varepsilon }(t),u)-l(t,\tilde{%
x}^{\varepsilon }(t),\tilde{y}^{\varepsilon }(t),\tilde{u}^{\varepsilon }(t))
\\
& -l(t,x^{\varepsilon }(t),y^{\varepsilon }(t),u)+l(t,x^{\varepsilon
}(t),y^{\varepsilon }(t),u^{\varepsilon }(t)) \\
& +\frac{1}{2}D^{2}(t)(u-\tilde{u}^{\varepsilon }(t))^{2}\tilde{P}%
_{1}^{\varepsilon }(t)-\frac{1}{2}D^{2}(t)(u-u^{\varepsilon
}(t))^{2}P_{1}^{\varepsilon }(t)\Big)\mbox{\rm d}t \\
\leq & \;C\varepsilon ^{\frac{\beta }{3}}.
\end{align*}%
Therefore, we get the first result on bounded control domains%
\begin{eqnarray*}
&&\int_{0}^{1}\left\langle p^{\varepsilon }(t),B(t)(u-u^{\varepsilon
}(t))\right\rangle +\left\langle k^{\varepsilon }(t),D(t)(u-u^{\varepsilon
}(t))\right\rangle -\left\langle q^{\varepsilon }(t),c(t)(u-\tilde{u}%
^{\varepsilon }(t))\right\rangle \\
&&+\frac{1}{2}D^{2}(t)(u-u^{\varepsilon }(t))^{2}P_{1}^{\varepsilon
}(t)+\theta _{0}^{\varepsilon }\left[ l(t,x^{\varepsilon }(t),y^{\varepsilon
}(t),u)-l(t,x^{\varepsilon }(t),y^{\varepsilon }(t),u^{\varepsilon }(t))%
\right] \mathrm{d}t \\
&\geq &-C\varepsilon ^{\beta }\theta _{0}^{\varepsilon },\text{ }\forall
u\in \mathbb{U},\text{ a.e, a.s.}.
\end{eqnarray*}

\noindent \textbf{Step 2. (The general case of control domains).}

For every $K=1,2,\cdots ,$ set
\begin{eqnarray*}
\mathbb{M}^{K} &\triangleq &\left\{ \left. y_{0}\in \mathbb{R}\right\vert
\left\vert y_{0}\right\vert \leq \left\vert y_{0}^{\varepsilon }\right\vert
+K\right\} , \\
\mathbb{N}^{K} &\triangleq &\left\{ \left. z\left( t\right) \in \mathbb{R}%
\right\vert \left\vert z\left( t\right) \right\vert \leq \left\vert
z^{\varepsilon }\left( t\right) \right\vert +K\right\} , \\
\mathcal{M}^{2}(0,1;\mathbb{N}^{K}) &\triangleq &\left\{ \left. z\left(
\cdot \right) \in \mathcal{M}^{2}(0,1;\mathbb{R})\right\vert z\left(
t\right) \in \mathbb{N}^{K}\right\} .
\end{eqnarray*}

\noindent Clearly, $\mathbb{M}^{K}$ is convex and $y_{0}^{\varepsilon }\in
\mathbb{M}^{K}\subseteq \mathbb{M}^{K+1},$ $\mathbb{R}=\cup _{K=1}^{\infty }%
\mathbb{M}^{K}.$ $z^{\varepsilon }(\cdot )\in \mathcal{M}^{2}(0,1;\mathbb{N}%
^{K})\subseteq \mathcal{M}^{2}(0,1;\mathbb{N}^{K+1}),$ and $\mathcal{M}%
^{2}(0,1;\mathbb{R})=\cup _{K=1}^{\infty }\mathcal{M}^{2}(0,1;\mathbb{N}%
^{K}).$ Note that $(y_{0}^{\varepsilon },z^{\varepsilon }(\cdot
),u^{\varepsilon }(\cdot ))$ is still a near optimal 3-triple of Problem $(%
\tilde{C}^{\varepsilon })$ when the original admissible control set is
replaced by $\mathbb{M}^{K}\times \mathcal{M}^{2}(0,1;\mathbb{N}^{K})\times
\mathcal{U}_{ad}\left[ 0,1\right] ,$ $K=1,2,\cdots .$ Moreover, (\ref{6.18})
also holds for fixed $\varepsilon >0$ on $\mathbb{M}^{K}\times \mathcal{M}%
^{2}(0,1;\mathbb{N}^{K})\times \mathcal{U}_{ad}\left[ 0,1\right] $ for every
$K=1,2,\cdots .$ Then there exists a subsequence
\begin{equation*}
\left( \theta _{0}^{\varepsilon ,K},\theta _{1}^{\varepsilon
,K},p^{\varepsilon ,K}(\cdot ),q^{\varepsilon ,K}(\cdot ),k^{\varepsilon
,K}(\cdot ),P^{\varepsilon ,K}(\cdot ),Q^{\varepsilon ,K}\left( \cdot
\right) \right)
\end{equation*}%
satisfying $\left\vert \theta _{0}^{\varepsilon ,K}\right\vert ^{2}+\mathbb{E%
}\left\vert \theta _{1}^{\varepsilon ,K}\right\vert ^{2}=1$, $\theta
_{0}^{\varepsilon ,K}\geq 0,$ (\ref{6.531})-(\ref{6.532}) such that the
following%
\begin{eqnarray*}
&&\int_{0}^{1}\left\langle p^{\varepsilon ,K}(t),B(t)(u-u^{\varepsilon
}(t))\right\rangle +\left\langle k^{\varepsilon ,K}(t),D(t)(u-u^{\varepsilon
}(t))\right\rangle \\
&&-\left\langle q^{\varepsilon ,K}(t),c(t)(u-\tilde{u}^{\varepsilon
}(t))\right\rangle +\frac{1}{2}D^{2}(t)(u-u^{\varepsilon
}(t))^{2}P_{1}^{\varepsilon ,K}(t) \\
&&+\theta _{0}^{\varepsilon ,K}\left[ l(t,x^{\varepsilon
,K}(t),y^{\varepsilon ,K}(t),u)-l(t,x^{\varepsilon ,K}(t),y^{\varepsilon
,K}(t),u^{\varepsilon }(t))\right] \mathrm{d}t\geq -C\varepsilon ^{\beta
}\theta _{0}^{\varepsilon ,K},\text{ }
\end{eqnarray*}%
holds. Since $\left\vert \theta _{0}^{\varepsilon ,K}\right\vert ^{2}+%
\mathbb{E}\left\vert \theta _{1}^{\varepsilon ,K}\right\vert ^{2}=1,$ there
is a subsequence also denoted by $\left( \theta _{0}^{\varepsilon ,K},\theta
_{1}^{\varepsilon ,K}\right) ,$ such that $\left( \theta _{0}^{\varepsilon
,K},\theta _{1}^{\varepsilon ,K}\right) \rightarrow \left( \theta
_{0}^{\varepsilon },\theta _{1}^{\varepsilon }\right) ,$ weakly in\textsl{\ }%
$\mathbb{R}\times L_{\mathcal{F}_{1}}^{2}(\Omega ;\mathbb{R})$, $\theta
_{0}^{\varepsilon }\geq 0.$ Hence, from continuous dependence of the
solution of BSDEs on parameters (see Yong and Zhou \cite{YZ}), we have
\begin{equation*}
\left( p^{\varepsilon ,K}(\cdot ),q^{\varepsilon ,K}(\cdot ),k^{\varepsilon
,K}(\cdot ),P^{\varepsilon ,K}(\cdot ),Q^{\varepsilon ,K}\left( \cdot
\right) \right) \rightarrow \left( p^{\varepsilon }(\cdot ),q^{\varepsilon
}(\cdot ),k^{\varepsilon }(\cdot ),P^{\varepsilon }(\cdot ),Q^{\varepsilon
}\left( \cdot \right) \right)
\end{equation*}%
in $\mathcal{M}^{2}(0,1;\mathbb{R})$ as $K\rightarrow +\infty .$ Moreover, $%
\left( p^{\varepsilon }(\cdot ),q^{\varepsilon }(\cdot ),k^{\varepsilon
}(\cdot ),P^{\varepsilon }(\cdot ),Q^{\varepsilon }\left( \cdot \right)
\right) $ satisfies (\ref{6.531})-(\ref{6.532}). Consequently, we get (\ref%
{3.3}). The proof is complete. ~\hfill $\Box $

\bigskip

\textbf{Acknowledgments.} \textsl{The authors would like to thank three
anonymous referees, AE, Professor Zhen Wu, Dr. Guangchen Wang, Dr. Zhiyong Yu for
their valuable comments, which led to a much better version of this article.
This article was partially done while the first author was visiting The Hong
Kong Polytechnic University in the summer of 2012. }

\end{document}